%% file: part2.tex
\documentclass[11pt]{amsart}
\usepackage{geometry}                
\usepackage{graphicx}
\usepackage{calc}
\usepackage{amssymb}
\usepackage{epstopdf}
\usepackage{hyperref}

\newcommand\rev[1]{#1}

\include{Var_Def}

\title{On Bayesian estimation and proximity operators}
\author{R. Gribonval and M. Nikolova}
\thanks{This work and the companion paper \cite{RGMN2018a} are dedicated to the memory of Mila Nikolova, who passed away prematurely in June 2018.
Mila dedicated much of her energy to bring the technical content to completion during the spring of 2018. The first author did his best to finalize the papers as Mila would have wished. He should be held responsible for any possible imperfection in the final manuscript.\protect\\
R. Gribonval, Univ Rennes, Inria, CNRS, IRISA, remi.gribonval@inria.fr; \protect\\
M. Nikolova, CMLA, CNRS and Ecole Normale Sup{\'e}rieure de Cachan, Universit{\'e} Paris-Saclay, 94235 Cachan, France.\\
\textcolor{blue}{Compared to the published version, this document (March 2026) includes typo corrections in Proposition~5, indicated in blue.}}

\begin{document}

\begin{abstract}
There are two major routes to address the ubiquitous family of inverse problems appearing in signal and image processing, such as denoising or deblurring. A first route relies on Bayesian modeling, where prior probabilities are used to embody models of both the distribution of the unknown variables and their statistical dependence with \rev{respect to} the observed data. The estimation process typically relies on the minimization of an expected loss (e.g. minimum mean squared error, or MMSE). The second route has received much attention in the context of sparse regularization and compressive sensing: it consists in designing (often convex) optimization problems involving the sum of a data fidelity term and a penalty term promoting certain types of unknowns (e.g., sparsity, promoted through an $\ell^1$ norm). 

Well known relations between these two approaches have \rev{led} to some widely spread misconceptions. In particular, while the so-called Maximum A Posterori (MAP) estimate with a Gaussian noise model does lead to an optimization problem with a quadratic data-fidelity term, we disprove through explicit examples the \rev{common} belief that the converse would be true. 

It has already been shown \cite{GRIBONVAL:2010:INRIA-00486840:1,NIPS2013_4868} that for denoising in the presence of additive Gaussian noise, for {\em any} prior probability on the unknowns, MMSE estimation can be expressed as a penalized least squares problem, with the apparent characteristics of a MAP estimation problem with Gaussian noise and a (generally) different prior on the unknowns. In other words, the variational approach is rich enough to build all possible MMSE estimators associated to additive Gaussian noise via a well chosen penalty.

We generalize these results beyond Gaussian denoising and characterize noise models for which the same phenomenon occurs. In particular, we prove that with (a variant of) {\em Poisson} noise and any prior probability on the unknowns, MMSE estimation can again be expressed as the solution of a penalized least squares optimization problem. 
For {\em additive} scalar denoising the phenomenon holds if and only if the noise distribution is log-concave. 
In particular, Laplacian denoising can (perhaps surprisingly) be expressed as the solution of a penalized least squares problem. 
In the multivariate case, the same phenomenon occurs when the noise model belongs to a particular subset of the exponential family. 
For multivariate {\em additive} denoising, the phenomenon holds if and only if the noise is white and Gaussian. 
\end{abstract}

\maketitle

\section{Introduction and overview}
\lab{sec:overview}

Inverse problems in signal and image processing consist in estimating an unknown signal $x_0$ given an indirect observation $y$ that may have suffered from blurring, noise, saturation, etc. The two main routes to address such problems are variational approaches and Bayesian estimation.

{\bf Variational approaches:} a signal estimate is the solution of an optimization problem
\beq \lab{var}
\hat{x} \in \arg\min_x D(x,y) + \ph(x)
\eeq
where $D(x,y)$ is a data-fidelity measure, and $\ph$ is a penalty promoting desirable properties of the estimate $\hat{x}$ such as, e.g., sparsity. 

A typical example is linear inverse problems, where one assumes $y = Lx_0 + e$ with $L$ some known linear operator (e.g., a blurring operator), $e$ is some error / noise. The most \rev{common} data-fidelity term is \rev{the square of the Euclidean norm $\|\cdot\|$}, which in combination with an $\ell^1$ sparsity-enforcing penalty leads to the well known Basis Pursuit Denoising approach
\beq \lab{varBP}
\hat{x}_{\text{BPDN}}(y) := \arg\min_x \tfrac{1}{2} \|y-Lx\|^2+\lambda \|x\|_1.
\eeq

{\bf Bayesian estimation:} $x_0$ is modeled as a realization of a random variable $X$ (with "prior" probability $p_X(x)$) and $y$ as the realization of a random variable $Y$ (with conditional probability distribution $p_{Y|X}(y|x)$). \rev{A Bayesian} estimator is designed as \rev{a} function $y \mapsto \hat{x}(y)$ that minimizes in expectation some specified cost $C(x_0,\hat{x})$, i.e., that minimizes
\beq \lab{bayesEst}
\mathbb{E}_{X,Y}\ C(X,\hat{x}(Y))
\eeq
where the pair $(X,Y)$ is drawn according to the joint distribution $p_{X,Y}(x,y) = p_{Y|X}(y|x) p_X(x)$. \rev{Equivalently, for a given $y$, the estimator $\hat{x}(y)$ is a minimizer of $\mathbb{E}_{X|Y=y} C(X,\hat{x}(y))$.}

A typical example is Minimum Mean Square Error (MMSE) estimation. The cost is the quadratic error $C(x,\hat{x}) = \|\hat{x}-x\|^2$, 
and the optimal estimator is the conditional expectation, also called Conditional Mean or Posterior Mean
\beq \lab{MMSE}
y \mapsto \hat{x}_{\text{MMSE}}(y) := \mathbb{E} (X|Y=y) = \int x p_{X|Y}(x|Y=y) dx.
\eeq
By Bayes law $p_{X|Y}(x|y) = \tfrac{p_X(x) p_{Y|X}(y|x)}{p_Y(y)}$, with $p_Y(y)$ the marginal distribution of $Y$.

While MMSE estimation yields the {\em expected value} of the {\em a posteriori} probability distribution $p_{X|Y}(x|y)$ of $X$, Maximum A Posteriori (MAP) estimation selects its mode, i.e. the {\em most probable} $x$ with respect to this distribution,
\[
\hat{x}_{\text{MAP}}(y) :=  
\arg\min_x -\log p_{X|Y}(x|y) = \arg\min_x \{ -\log p_{Y|X}(y|x)-\log p_X(x)\}.
\]
MAP is directly connected to variational approaches, hence its popularity. However, this is not usually considered as a proper Bayesian estimator although it can be seen as minimizing \eq{bayesEst} with the "pseudo-cost" $C(x,\hat{x}) := \delta(x-\hat{x})$. We will soon come back to the Bayesian interpretation of MAP estimation and its pitfalls.

Many other costs can be used, e.g. $C(x,\hat{x}) = \|x-\hat{x}\|$ yields the conditional spatial median.

Banerjee et al \cite{Banerjee:2005jd} show that if $C(x,\hat{x})$ (defined on $\RR^n \times \RR^n$) is the Bregman divergence\footnote{By definition $D_\hh(x,y) :=  \hh(x)-\hh(y) -\<\nabla\hh(y),x-y\>$. This is usually not symmetric, i.e. $D_{h}(x,y) \neq D_{h}(y,x)$.} $D_h(x,\hat{x})
$ of a strictly convex proper differentiable function $h$ \cite{Bregman67}, then the conditional mean is the unique minimizer of \eq{bayesEst} \cite[Theorem 1]{Banerjee:2005jd}. Vice-versa, they prove under mild smoothness assumptions on $C(x,\hat{x})$ that if the conditional mean is the unique minimizer of \eq{bayesEst} for any pair of random variables $X,Y$ then 
\beq \lab{banerjeeCM}
C(x,\hat{x}) = D_h(x,\hat{x}),\qu \all x,\hat{x}
\eeq 
for some strictly convex differentiable function $h$. 

\subsection{The MAP vs MMSE quid pro quo} 
A \rev{common} quid pro quo between tenants of the two approaches revolves around the MAP interpretation of variational approaches \cite{SSP}.  
In the particular case of a linear inverse problem with white Gaussian noise the conditional density reads\rev{\footnote{\rev{with $\propto$ denoting proportionality.}}} $p_{Y|X}(y|x) \propto \exp\lp(-\tfrac{\|y-Lx\|^2}{2\sigma^2}\rp)$, and denoting $\ph_X(x) := -\sigma^2 \log p_X(x)$, MAP estimation reads
\beq \lab{MAPisVAR}
\hat{x}_{\text{MAP}} = \arg\min_x \tfrac{1}{2\sigma^2} \|y-Lx\|^2 - \log p_X(x) = \arg\min_x \tfrac{1}{2} \|y-Lx\|^2 + \ph_X(x)
\eeq
Thus, {\em if} one assumes a Gaussian noise model and {\em if} one chooses MAP as an estimation principle {\em then} this results in a variational problem shaped as \eq{var} with a quadratic data-fidelity term and a penalty which is a scaled version of the negative log-prior. 

\rev{It is argued in} 
 \cite{Nikolova:2007aa} that --except in very special circumstances (Gaussian prior and Gaussian noise)-- ``{\em the MAP approach is not relevant in the applications where the data-observation and the prior models are accurate}'', \rev{in that}
 \eq{MAPisVAR} leads to a suboptimal estimator when the considered data is indeed distributed as $p_X(x) \propto \exp(-\ph_X(x)/\sigma^2)$ with Gaussian noise. 
 \rev{Indeed, consider as an example} compressive sensing where $y = Lx+e$ with Gaussian i.i.d. $e$ and $L$ an underdetermined measurement matrix. When $X$ is a random vector with i.i.d. Laplacian entries, we get $\ph_X(x) \propto \|x\|_1$ and the MAP estimator is Basis Pursuit Denoising \eq{varBP} which has been shown \cite{gribonval:inria-00563207} to have poorer performance (in the highly compressed regime and in the limit of low noise) than a variational estimator \eq{var} with quadratic data-fidelity and quadratic penalty $\ph(x) \propto \|x\|_2^2$, aka Tikhonov regression or ridge regression.

Unfortunately, a widely spread misconception has \rev{led} to a "reverse reading" of optimization problems associated to variational approaches. 
For example, even though it is true that one obtains \eq{varBP} as a MAP under additive Gaussian noise with a Laplacian signal prior, by no means does this imply that the use of \eq{varBP} to build an estimator is necessarily motivated by the {\em choice} of MAP as an estimation principle and the {\em belief} that the Laplacian prior is a good description of the distribution of $X$. 
Instead, as widely documented in the literature on sparse regularization, the main reason for choosing the $\ell^1$ penalty is simply to promote sparse solutions: 
any minimizer of \eq{varBP} is bound to have (many) zero entries (in particular when the parameter $\lambda$ is large) which is desirable when prior knowledge indicates that $x_0$ is ``compressible", that is to say, well approximated by a sparse vector. 

As demonstrated e.g. in \cite{gribonval:inria-00563207} (see also \cite{Amini:2011uh} for related results), a random vector $X \in \RR^n$ with entries drawn i.i.d. from a Laplacian distribution is, with high probability, {\em not} compressible; on the opposite, when the entries of $X$ are drawn i.i.d. according to a {\em heavy-tailed distribution}, $X$ is typically compressible. Thus, despite having the apparent characteristics of the MAP estimator with a Laplacian prior, \eq{varBP} in fact approximates well the MMSE estimator with a heavy-tailed prior.

\rev{\begin{remark}\label{rk:MAPorMMSEbetter}
One should be careful about the temptation of determining once and for all which of MAP or MMSE is a ``better'' estimation principle. While the above compressive sensing example illustrates that MAP estimation under the true underlying distribution can lead to poorer estimation performance than MMSE, J. Idier (private correspondence) provided the authors with an interesting  example of the converse phenomenon, where MMSE estimation is not as operational as MAP. This happens, e.g. in phase unwrapping where $x$ needs to be estimated modulo $2\pi$. MAP then involves finding one of the modes of the posterior $p_{X|Y}(x|y)$, all modes being more or less equal modulo $2\pi$, hence equally operational. On the contrary MMSE typically averages several modes and leads to non-operational estimates. Similar problems may arise in other estimation problems where label switching may occur, leading to multiple modes of the posterior distribution. \\
\end{remark}
\begin{remark}\label{rk:GoodEstimatorEqualDist}
In  \cite{Nikolova:2007aa}, an argument against the MAP approach is expressed as:
{\em "In full rigor, an estimator $\hat{X}$ for $X$, based on data $Y$, can be said to be coherent with the underlying models if $\hat{X} \sim f_X$"}, i.e., if $\hat{X}$ has the same probability distribution as $X$. Again, as pointed out by J. Idier, such a criterion to qualify a ``good'' estimator seems questionable: while ``$\hat{X} \sim f_{X}$'' can be (approximately) expected when the knowledge of $Y$ allows to perfectly estimate $X$ (or, on the opposite, when $Y$ brings no information on $X$!), such a property cannot be expected in the more realistic intermediate cases: even in the linear Gaussian setting, where MAP and MMSE lead to the same (unbiased) estimator, the covariance of  the estimator generally does not match that of the prior. A more convincing intrinsic issue with MAP is probably its lack of invariance with respect to reparametrization of the problem.
\end{remark}}

\subsection{Writing certain {\em convex} variational estimators as proper Bayes estimators} 
In penalized least squares regression for linear inverse problems in $\RR^n$, the variational estimator
\beq \lab{penLS}
\hat{x}(y) \in \arg \min_{x \in \RR^n} \tfrac{1}{2} \|y-Lx\|^2 + \ph(x),
\eeq
is traditionally interpreted as a MAP estimator \eq{MAPisVAR} with Gaussian noise and \rev{prior}
\[
p_X(x) \propto \exp(-\ph(x)).
\]
When the penalty $\ph$ is convex differentiable and $\tfrac{1}{2}\|Lx\|^2+\ph(x)$ has at least linear growth at infinity, by \cite[Theorem 1]{Burger:2014hv} the estimator \eq{penLS} is also a proper Bayesian estimator minimizing \eq{bayesEst}, with a prior distribution $p_X(x) \propto \exp(-\ph(x))$ and $y = Lx+e$ where $e$ is white Gaussian noise, for the cost
\beq \lab{burgerCost}
C(x,\hat{x}) = \tfrac{1}{2}\|L(x-\hat{x})\|^2+D_{\ph}(\hat{x},x) 
\eeq
involving the Bregman divergence $D_\ph$. 
The results in \cite{Burger:2014hv} are actually expressed with colored Gaussian noise and \rev{with} $\|\cdot\|$ replaced by the corresponding weighted (Mahalanobis) quadratic norm in \eq{penLS} and \eq{burgerCost}. Further extensions to infinite-dimension have been considered in \cite{Helin:2015jq}.

In other words, in the context of linear inverse problems with additive Gaussian noise and a \rev{log-concave} prior, a variational approach \eq{penLS} having all the apparent features of the MAP estimator is in fact also a proper Bayesian estimator with the specific cost \eq{burgerCost}.

 \BR
 The reader may notice that the cost \eq{burgerCost} is a Bregman divergence $C(x,\hat{x}) = D_h(\hat{x},x)$, with $h(x) := \tfrac{1}{2} \|Lx\|^2+\ph(x)$.
 The order of arguments is reversed compared to \eq{banerjeeCM}, which is neither an error nor a coincidence: by the results of Banerjee \cite{Banerjee:2005jd}, when $C(x,\hat{x}) = D_h(x,\hat{x})$ for some $h$, the  Bayes estimator minimizing \eq{bayesEst} is simply the MMSE. Here, in \eq{penLS}-\eq{burgerCost}, 
 this does not happen: the MAP estimator does not match the MMSE estimator except in certain very specific cases where the Bregman divergence is indeed symmetric.
\ER

\subsection{Writing certain MMSE estimators using an expression analog to a MAP estimator} \lab{bayes}

The results of Burger et al strongly \rev{intertwine} the prior model $p_X(x) \propto e^{-\ph(x)}$, the observation model $p_{Y|X}(y|x)$ embodied by the linear operator $L$, and the task cost $C_{L,\ph}(x,\hat{x})$. In particular, the task ("what we want to do") becomes {\em dependent} on the data model ("what we believe"). 

From a data processing perspective the above approach is not fully satisfactory: it seems more natural to first choose a relevant task (e.g., MMSE estimation) and a reasonable model (e.g., a prior $p_X(x)$), and then to design a penalty (the tool used to solve the task given the model) based on these choices.

In this spirit, in the context of additive white Gaussian denoising, \cite{GRIBONVAL:2010:INRIA-00486840:1} showed that, {\em for any signal prior $p_X(x)$}, the MMSE estimator \eq{MMSE} is the unique solution (and unique stationary point) of a variational optimization problem 
\beq \lab{pseudoMAP}
\hat{x}_{\text{MMSE}}(y) = \arg\min_x \tfrac{1}{2}\|y-x\|^2+\wt{\ph}_X(x)
\eeq
with $\wt{\ph}_X(x)$ some penalty that depends on the considered signal prior $p_X(x)$. 

In other words, MMSE has all the apparent features of a MAP estimator with Gaussian noise and a "pseudo" signal prior $\wt{p}_X(x) \propto e^{-\wt{\ph}_X(x)}$. Except in the very special case of a Gaussian prior,  {\em the pseudo-prior $\wt{p}_X(x)$ differs from the prior $p_X(x)$} that defines the MMSE estimator \eq{MMSE}. This result has been extended to MMSE estimation for inverse problems with additive colored Gaussian noise \cite{NIPS2013_4868}. Unser and co-authors \cite{Kazerouni:2013co,Amini:2013co},\rev{\cite[Section 10.4.3]{Unser:2014vs}} exploit these results for certain MMSE estimation problems. Louchet and Moisan \cite[Theorem 3.9]{Louchet:2013hs} consider the specific case of the MMSE estimator associated to a total variation image prior $p_X(x)$ and establish the same property through connections with the notion of {\em proximity operator} of a convex lsc function. 

\subsection{Contribution: MMSE estimators that can be expressed as proximity operators}
We extend the general results of \cite{gribonval:inria-00563207,NIPS2013_4868} beyond Gaussian denoising using a characterization of proximity operators of possibly nonconvex penalties obtained in a companion paper \cite{RGMN2018a}.
Our extension goes substantially beyond Gaussian denoising, including scalar Poisson denoising, scalar denoising in the presence of additive noise with any log-concave distribution, and multivariate denoising for certain noise distributions belonging to the exponential family.

\subsubsection{Scalar denoising}
\BP[scalar Poisson denoising]\lab{scalpoisson}
Consider the scalar Poisson noise model where the conditional probability distribution of the integer random variable $Y \in \NN$ given $x \in \RR_+^*$ is
\[
p_{Y|X}(Y=n|x) = \tfrac{x^n}{n!} e^{-x},\qu \all n \in \NN
\]
and let $p_X$ be any probability distribution for a positive random variable $X>0$. 
There is a (possibly nonconvex)  penalty function $\tilde{\ph}_X: \RR_+^* \to \RR \cup \{+\infty\}$ such that, for any $n \in \NN$
\beq \lab{scalpoissonprox}
\EXP(X|Y=n) \rev{~\in~} \arg\min_{x>0} \lp\{\tfrac{1}{2}(n-x)^2+\tilde{\ph}_X(x)\rp\}.
\eeq
Note that the $\arg\min$ is strictly positive since the conditional expectation is strictly positive.
The penalty $\tilde{\ph}_X$ depends \rev{on the probability distribution $p_X$ of $X$.}
\EP

\BP[scalar additive noise]\lab{scallogconc}
Consider an {\em additive} noise model $Y=X+N$ where the random variables $X,Z \in \RR$ are independent. The conditional probability distribution of the random variable $Y \in \RR$ given $x \in \RR$ is $p_{Y|X}(y|x) = p_Z(y-x)$. Assume that $p_Z(z)>0$ for any $z \in \RR$ and that $z \mapsto F(z) := -\log p_Z(z)$ is continuous. The following properties are equivalent:
\bena
\item \lab{scallca} the function $F$ is convex (i.e., the noise distribution is log-concave);
\item \lab{scallcb} for any prior probability distribution $p_X$ on the random variable $X \in \RR$,  the conditional expectation $\EXP(X|Y=y)$ is well defined for any $y \in \RR$, and there is a (possibly nonconvex)  penalty function $\tilde{\ph}_X: \RR_+^* \to \RR \cup \{+\infty\}$ such that for any $y \in \RR$
\beq \lab{scallogconcprox}
\EXP(X|Y=y) \rev{~\in~} \arg\min_{x \in \RR} \lp\{\tfrac{1}{2}(y-x)^2+\tilde{\ph}_X(x)\rp\}.
\eeq
The penalty $\tilde{\ph}_X$ depends \rev{on the probability distribution $p_X$ of $X$ and the function $F$.}
\een
\EP

\prop{scalpoisson} and \prop{scallogconc} are proved in Section~\ref{mmse} as corollaries of a more general result (\lem{ScalMMSE}) on scalar MMSE estimation.
Let us insist that while the naive interpretation of an optimization problem such as \eq{scalpoissonprox} (resp.~\eq{scallogconcprox}) would be that of MAP estimation with Gaussian noise, it actually corresponds to MMSE estimation with Poisson (resp. log-concave) noise.

Classical log-concave examples include generalized Gaussian noise \cite{Mathieu92,Belge00}, where $p_{Y|X}(y|x) \propto \exp\lp(-\lp|\tfrac{x-y}{\sigma}\rp|^\gamma\rp)$ with $1 \leq \gamma < \infty$. This includes Gaussian noise for $\gamma=2$, but also Laplacian noise for $\gamma=1$, and in all cases the MMSE estimator can always be written as a "pseudo-MAP" with a "Gaussian" (i.e., quadratic) data-fidelity term and an adapted "pseudo-prior" $\tilde{\ph}_X$. For $0<\gamma <1$ the noise distribution is not log-concave, hence \rev{\em there are prior probability distributions $p_{X}$ such that} the corresponding MMSE estimator \rev{\em cannot} be written as in \eq{scallogconcprox}.

\begin{figure}[htbp]
\includegraphics[width=\textwidth*2/3]{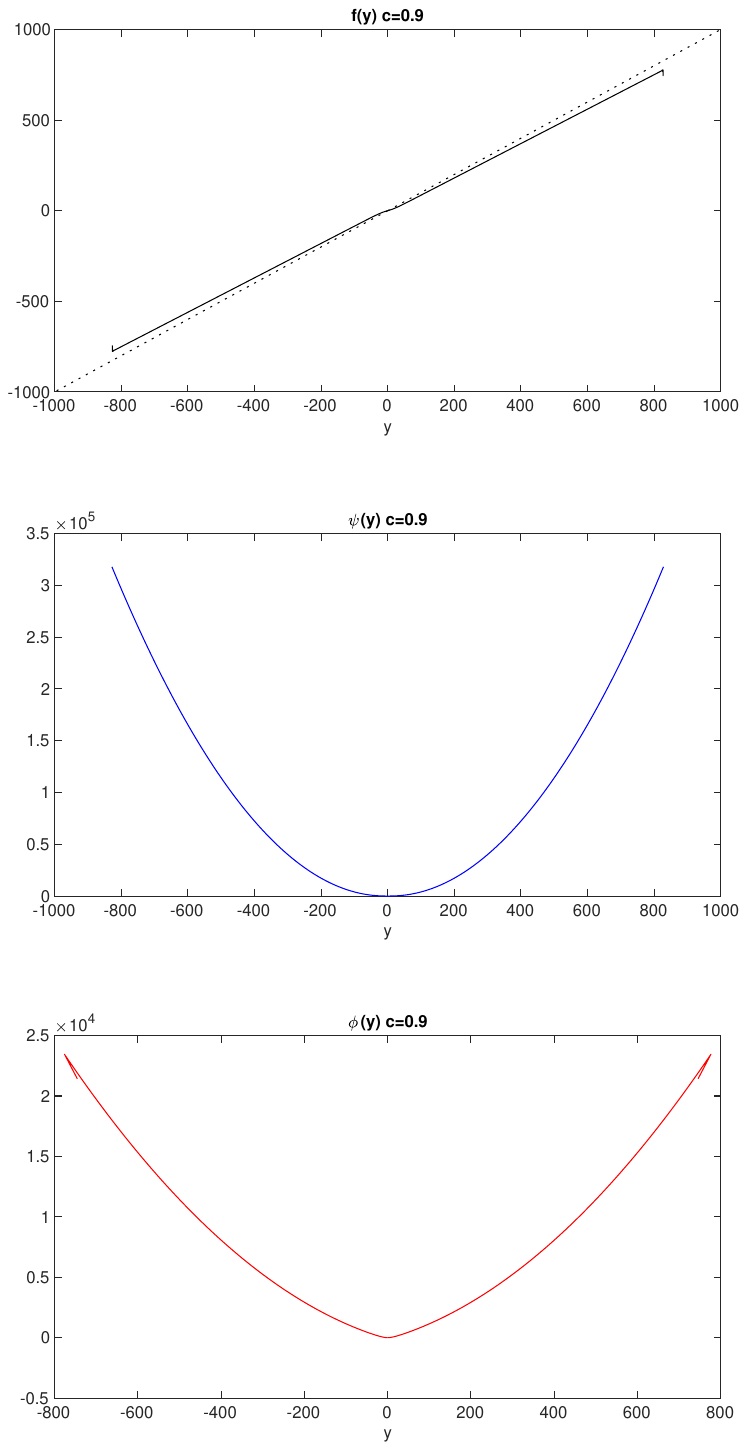}
\caption{$f(y)$ (top), potential $\psi(y)$ \rev{from Theorem~\ref{ProCon}} (middle) and penalty $\ph(y)$ (bottom) for \exa{scalL1L1} with $c=0.9$.}
\label{figL1L1} 
\end{figure}

\BE \lab{scalL1L1}
Consider $X$ a scalar random variable with a Laplacian prior distribution, and $Y=X+cZ$ where $Z$ is an independent scalar Laplacian noise variable and $c>0$. While the MAP estimator reads
\[
\arg\min_x \lp\{|y-x|+c|x|\rp\},
\]
by \prop{scallogconc} the MMSE estimator can be expressed as
\[
f(y) = \arg\min_x \lp\{\tfrac{1}{2} (y-x)^2+\ph(x)\rp\}
\]
with some penalty $\ph$.  The details of the analytic derivation of $f$ are in Appendix \ref{appscalL1L1}. The corresponding potential $\psi$ \rev{(cf Theorem~\ref{ProCon} in Section~\ref{mmse})} can be obtained by numerical integration, and the penalty $\ph$ is characterized by  \rev{Eq.~}\eq{psivsphi}: $\ph(f(y))=yf(y)-f^2(y)/2-\psi(y)$. It can be plotted using the pairs $f(y)$,$\ph(f(y))$. Figure~\ref{figL1L1} provides the shape of $f(y)$, of the corresponding potential $\psi(y)$, and of the corresponding penalty $\ph(x)$ for $c=0.9$. 
\EE

\subsubsection{Multivariate denoising}

 The above scalar examples straightforwardly extend to the multivariate setting in the special case where both the noise model and the prior $p_X$ are completely separable, cf \exa{compsep} in Section~\ref{mmse}. 
In general however, the multivariate setting is quite different from the scalar one.
For example in $\RR^n$, $n \geq 2$, MMSE estimation in the presence of additive Laplacian noise (resp. of Poisson noise) {\em cannot always} be written as in \eq{pseudoMAP}, depending on the prior $p_X$, see \exa{exlapv} (resp. \exa{poisson}) in Section~\ref{mmse}.
In contrast, we show (\lem{VectExpFam}) that the MMSE estimator can always be expressed as in \eq{pseudoMAP}, \rev{provided we consider} particular noise models of the exponential family
\beq \lab{expfam0}
p_{Y|X}(y|x) = \exp\lp(c\<x,y\>-a(x)-b(y)\rp)
\eeq
for some $c \geq 0$ and smooth enough $b$. This form is in fact essentially necessary (\lem{VectMMSENec}).

A primary example is additive Gaussian white noise, where 
\beq \lab{gaussnoise}
p_{Y|X}(y|x) = C \exp\lp(-\tfrac{\|y-x\|^2}{2\sigma^2}\rp) = \exp\lp(\log C -\tfrac{\|y\|^2}{2\sigma^2}-\tfrac{\|x\|^2}{2\sigma^2}+\tfrac{1}{\sigma^2}\<x,y\>\rp),
\eeq
and we recover the main results of \cite{GRIBONVAL:2010:INRIA-00486840:1}, using the following definition.
\BD \lab{nondeg}
A random variable $X$ with values in \rev{a Hilbert space} $\H$ \rev{and probability distribution $P$} is non-degenerate if there is no affine hyperplane $\V \subset \H$ such that $X \in \V$ almost surely, i.e., if $\rev{P}(\<X,v\>=d)<1$ for all nonzero $v \in \H$ and $d \in \RR$.
\ED

\BP[additive white Gaussian denoising] \lab{gauss}
Consider an additive noise model $Y=X+Z \in \RR^n$ where the random variables $X,Z$ are independent and $p_Z(z) 
\propto \exp\lp(-\tfrac{\|z\|^2}{2\sigma^2}\rp)$. Then for any prior distribution $p_X$ there is a (possibly nonconvex) penalty $\wt{\ph}_X: \RR^n \to \RR \cup \{+\infty\}$ such that for any $y \in \RR^n$
\beq \lab{gaussprox}
\EXP (X|Y=y) \rev{~\in~} \arg\min_{x \in \RR^n} \{\tfrac{1}{2} \|y-x\|^2+\wt{\ph}_X(x)\}.
\eeq
The penalty $\wt{\ph}_X$ depends \rev{on the probability distribution $p_X$ of $X$ and the variance $\sigma^{2}$}. 

Further, if $X$ is non-degenerate then the function $y \mapsto \EXP(X|Y=y)$ is injective and there is a choice of $\wt{\ph}_X$ such that for any $y$, $\EXP(X|Y=y)$ is the unique stationary point (and global minimizer) of the rhs of \eq{gaussprox}, \rev{hence $y \mapsto \EXP(X|Y=y)$ {\em is} the proximity operator of $\wt{\ph}_{X}$.}
\EP
{\remark It is further known \cite{GRIBONVAL:2010:INRIA-00486840:1} that in this Gaussian context the conditional mean estimator \rev{$y \mapsto \EXP(X|Y=y)$} is $C^\infty$, and that $e^{-\wt{\ph}_X(x)}$ is (up to renormalization) a proper prior density. It is also known \cite{NIPS2013_4868} that $\wt{\ph}_X$ is convex (resp. additively separable) if and only if the marginal $p_Y(y)$ is log-concave (resp. separable); this holds in particular as soon as the prior $p_X$ is log-concave (resp. separable). Extending such characterizations beyond Gaussian denoising is postponed to further work.
}

This is the only example of centered additive noise with smooth density of the form \eq{expfam0}.  
\BP \lab{vaddmv}
Consider a multivariate {\em centered additive} noise model i.e. with $p_{Y|X}(y|x) = p_Z(y-x)$. Assume that $p_Z(z)>0$ for any $z \in \H$ and denote $z \mapsto F(z) := - \log p_Z(z)$. If $F$ is continuous and $p_{Y|X}$ is of the form \eq{expfam0} then $p_{Y|X}(y|x)$ is in fact Gaussian, with the expression \eq{gaussnoise} for some $\sigma>0$.
\EP
The proof is in Appendix~\ref{pfvaddmv}. Despite \prop{vaddmv} and the apparent connection between the quadratic fidelity term in \eq{gaussprox} and the Gaussian log-likelihood $-\log p_{Y|X}(y|x)$, noise models of the form \eq{expfam0} are more general. For example, they cover a variant of multivariate Poisson denoising.

\BP[variant of multivariate Poisson denoising]\lab{logpoisson}
Consider a random variable $Y \in \NN^n$ with conditional distribution given $x \in (\RR_+^*)^n$ expressed as
\[
p_{Y|X}(Y=y|x) = \prod_{i=1}^n \tfrac{x_i^{y_i}}{y_i!} e^{-x_i},\qu \all y \in \NN^n,
\]
and let $p_X$ be any probability distribution for a multivariate positive random variable $X$. 
There is a (possibly nonconvex)  penalty function $\tilde{\ph}_X: \RR^n \to \RR \cup \{+\infty\}$ such that the MMSE estimator {\em of the (entrywise) logarithm} $\log X = (\log X_i)_{i=1}^n$ satisfies for any $y \in \NN^n$
\beq \lab{logpoissonprox}
\EXP(\log X|Y=y) \rev{~\in~} \arg\min_{\xi \in \RR^n} \lp\{\tfrac{1}{2}\|y-\xi\|^2+\tilde{\ph}_X(\xi)\rp\}.
\eeq
The penalty $\tilde{\ph}_X$ depends \rev{on the probability distribution $p_X$ of $X$}. 

Further, if $\textcolor{blue}{\log}\, X$ is non-degenerate then the function $y \mapsto \EXP(\textcolor{blue}{\log}\, X|Y=y)$ is injective and there is a choice of $\wt{\ph}_X$ such that for any $y$, $\EXP(\textcolor{blue}{\log}\, X|Y=y)$ is the unique stationary point (and global minimizer) of the rhs of \eq{logpoissonprox},  \rev{hence $y \mapsto \EXP(\textcolor{blue}{\log}\, X|Y=y)$ {\em is} the proximity operator of $\wt{\ph}_{X}$.}
\EP
In \prop{logpoisson} the noise model is not Gaussian, yet MMSE estimation can be expressed in a variational form "looking like" Gaussian denoising due to the quadratic data-fidelity term in \eq{logpoissonprox}. 
\prop{gauss} and \prop{logpoisson} are proved in Section~\ref{mmse} as corollaries of \lem{VectExpFam}.

\subsection{Discussion and extensions}
While this paper primarily focuses on characterizing when the conditional mean estimator is a proximity operator, it is natural to wonder under which noise models the conditional mean estimator can be expressed as the solution of another variational optimization problem such as \eq{varBP} or more generally \eq{var}, with properly chosen data fidelity term and penalty. This has been done \cite{NIPS2013_4868} for linear inverse problems with colored Gaussian noise, and we expect that it is possible to combine the corresponding proof techniques with \cite[Theorem 3(c) \rev{and Corollary 6}]{RGMN2018a} to obtain results in this vein for a larger class of observation models. Of particular interest would be to understand whether Poisson denoising can be written as the solution of a variational problem \eq{var} with a Kullback-Leibler divergence as data-fidelity term (which appears naturally in a MAP framework) and a well chosen penalty.

Finally, the reader may have noticed that the characterizations obtained in this paper are non constructive. They merely state the {\em existence} of a penalty such that the MMSE, which is a priori expressed as an integral, is in fact the solution of a (often convex) variational problem. From a practical point of view, a challenging perspective would be to identify how to exploit this property to design efficient estimation algorithms.

\rev{\paragraph{\bf Acknowledgement}
The corresponding author is indebted to the anonymous reviewers for their remarks, as well as to J. Idier for his insightful comments on a preliminary version of this work which have led in particular to Remarks~\ref{rk:MAPorMMSEbetter} and \ref{rk:GoodEstimatorEqualDist}.}

\section{When is MMSE estimation a proximity operator ?}
\label{mmse}

We now provide our main general results on the connections between Bayesian estimation and variational approaches.
After some reminders on the expression of the conditional mean in terms of a marginal and a conditional distribution, we define a class of "proto" conditional means estimators based on a given "proto" conditional distribution $\qyx{y}{x}$. Then, we focus on the scalar case where we characterize (proto) conditional distributions that lead to (proto) conditional mean estimators that are proximity operators. Finally we consider the vector case.

\subsection{Reminders on proximity operators}
\rev{ \cite{Helin:2015jq}}

Let $\H$ be a Hilbert space equipped with an inner product denoted $\<\cdot,\cdot\>$ and a norm denoted $\|\cdot\|$. In the finite dimensional case, one may simply consider $\H = \RR^{n}$. A function $\ff : \H \to \RR \cup \{+\infty\}$ is proper iff there is $x\in\H$ such that $\ff(x) < +\infty$, i.e., $\dom{\ff} \neq \void$, where $\dom{\ff} := \{x \in \H\mid \ff(x) < \infty\}$. It is lower semi-continuous (lsc) if for any $x_0 \in \H$, $\liminf_{x \to x_0} \ff(x) \geq \ff(x_0)$, or equivalently if the set $\{x \in \H: \ff(x)>\alpha\}$ is open for every $\alpha \in \RR$. 
The proximity operator of a (possibly nonconvex) proper penalty function $\ph: \H \to \RR \cup \{+\infty\}$ is the set-valued operator
\[
y \mapsto \prox_{\ph}(y) := \arg\min_{x \rev{\in\H}} \lp\{ \tfrac12 \|y-x\|^{2}+\ph(x) \rp\}
\]
A primary example is soft-thresholding $f(y) := y(1-1/y)_+$, $y \in \H:=\RR$ , which is the proximity operator of the absolute value function $\ph(x):=|x|$. Let us recall that if $\ff$ is a convex function, $u \in \H$ is a subgradient of $\ff$ at $x \in \H$ iff $\ff(x')-\ff(x) \geq \<u,x'-x\>$, $\all x' \in \H$. The subdifferential at $x$ is the collection of all subgradients of $\ff$ at $x$ and is denoted $\d \ff(x)$.

The following results are proved in the companion paper \cite{RGMN2018a}. \rev{Strictly speaking, a proximity operator is set-valued: a function $f$ such that $f(y) \in \prox_\ph(y)$ for any $y \in \Y$ is a {\em selection} of the proximity operator of $\varphi$. For concision we will say that $f$ {\em implements} this proximity operator.}
\BT\cite[Theorem 1]{RGMN2018a}\lab{ProCon}
Consider $\Y$ a non-empty subset of $\H$.
A function $f: \Y\to \H$ \rev{implements} the proximity operator of some penalty $\ph$ (i.e.
$f(y) \in \prox_\ph(y)$ for any $y \in \Y$) iff there exists a convex lsc function $\psi: \H \to \RR \cup \{+\infty\}$ such that for any $y \in \Y$, $f(y) \in \d \psi(y)$. 
\ET
When the domain $\Y$ is convex, \cite[Theorem 3]{RGMN2018a} implies that  there is a number $K \in \RR$ such that the functions $f$, $\ph$ and $\psi$ in \thm{ProCon} satisfy 
\beq
\lab{psivsphi}
\psi(y) = \<y,f(y)\> -\tfrac{1}{2}\|f(y)\|^2-\ph(f(y)) + K,\qu \all y \in \Y.
\eeq
\BC\cite[Corollary 4]{RGMN2018a} \lab{coscal}
Consider an arbitrary non-empty subset $\Y \subset \RR$. A function $f: \Y \to \RR$ \rev{implements} the proximity operator of some penalty $\ph$ if, and only if, it is  nondecreasing. 
\EC
\BT\cite[Theorem 2]{RGMN2018a}\label{jac1}
Let $\Y$ be an open convex subset of $\H$ and $f : \Y \to \H$ be $C^1(\Y)$.
The following properties are equivalent:
\bena
\item\label{jac0} $f$ \rev{implements} the proximity operator of a function $\ph$ (i.e. $f(y) \in \prox_{\ph}(y)$ for any $y\in\Y$);
\item\label{jaca} there is a convex $C^{2}$ function $\psi: \H \to \RR \cup \{+\infty\}$ such that $f(y) = \nabla \psi(y)$ for all $y \in \Y$;
\item \label{jacb} the differential $Df(y)$ is a {\em symmetric positive semi-definite} operator\footnote{\label{foot1}
A continuous linear operator $L: \H \to \H$ is symmetric if $\<x, Ly\> = \<Lx, y\>$ for any $x,y \in \H$. A symmetric continuous linear operator is positive semi-definite if $\<x,Lx\> \geq 0$ for any $x \in \H$. This is denoted $L \succeq 0$. It is positive definite if $\<x,Lx\> >0$ for any nonzero $x \in \H$. This is denoted $L \succ 0$.} 
for any $y \in \Y$.
\een
\ET
\BC\cite[Corollary 3]{RGMN2018a}\lab{uniqueC1}
Let $\Y \subset \H$ be open and convex, and $f: \Y \to \H$ be $C^1$ with $Df(y) \succ 0$ for any $y \in \Y$. 
Then $f$ is injective and there is $\ph: \H \to \RR \cup \{+\infty\}$ such that $\prox_\ph(y) = \{f(y)\},\ \all y \in \Y$ and $\dom \ph = \Ima{f}$. Moreover, if $x \in \H$ is a stationary point\footnote{$u$ is a stationary point of $\gg: \H \to \RR \cup \{+\infty\}$ if $\nabla \gg(u)=0$; then, $\gg$ is proper on a neighborhood of $u$.} of $x \mapsto \tfrac{1}{2}\|y-x\|^2+\ph(x)$
then $x=f(y)$.
\EC

\subsection{Reminders on conditional expectation}
Consider a pair of random variables $(X,Y)$ with values in
\rev{$\RR^{n} \times \RR^{m}$}, with joint probability density function (pdf)
$p_{X,Y}(x,y)$ and marginals $p_{Y}(y),p_{X}(x)$. 
For $y$ such that $p_{Y}(y)>0$, the conditional distribution of $X$ given $Y=y$ is $p_{X|Y}(x|y) = p_{X,Y}(x,y)/p_{Y}(y)$. When $\|x\|_{2}$ is integrable with respect to $p_{X|Y}(\cdot|y)$, the conditional expectation of $X$ given $Y=y$ is
\[\EXP(X|Y=y) = \int x\ p_{X|Y}(x|y)dx\]
By Bayes rule, the conditional distribution and the marginal $p_Y(y)$ satisfy
\[p_{X|Y}(x|y) = \frac{p_{Y|X}(y|x) p_{X}(x)} {p_{Y}(y)} \qu \mb{and} \qu p_{Y}(y) = \int p_{X,Y}(x,y)dxdy 
= \int p_{Y|X}(y|x) p_{X}(x) dx.
\]
Denoting
\[\qyx{y}{x} := p_{Y|X}(y|x)\]  
we thus have 
\beqnn
\EXP(X|Y=y) &=& \frac{\int x\ p_{Y|X}(y|x) p_X(x) dx}{p_Y(y)} = 
\frac{\EXP_X\ 
\big(X \qyx{y}{X}\big)}{\EXP_X\ \big( \qyx{y}{X}\big)}.
\eeqnn
As we will see (\cor{scalpoisson}), the conditional mean has the same expression in related settings such as scalar Poisson denoising. Considering a function $\qyx{y}{x}$ that plays the role of a "proto" conditional distribution of the observation $y$ given the unknown $x$, we can define "proto" conditional expectation functions in order to characterize when the conditional expectation \rev{implements} a proximity operator. 
\subsection{"Proto" conditional distributions and "proto" conditional expectations}

\BD\lab{def:ProtoMMSE}
Consider \rev{a Hilbert space $\H$,} $\X,\Y \subset \H$ and a "proto" conditional distribution $q : \X \times \Y \to \RR_+, (x,y) \mapsto \qyx{y}{x}$. Given \rev{a random variable $X$ with distribution $P$} on $\X$ such that
\beq\label{eq:admrv}
\EXP_{X \sim P}\ \big((1+\|X\|) \qyx{y}{X}\big) < \infty,\qu  \all y \in \Y. 
\eeq
The "proto" marginal distribution $q_P(y)$ of $\qyx{y}{x}$ is defined\footnote{If \rev{$\H = \RR^{n}$ and } $p_{Y|X}(y|x) = \qyx{y}{x}$ is a well-defined conditional probability then $q_P(y) = \int p_{X,Y}(x,y) dx$ is simply the marginal distribution of the random variable $Y$ where $p_{X,Y}(x,y) = p_{Y|X}(y|x) P(x) = \qyx{y}{x} P(x)$.} as the function 
\[
y\in\Y \mapsto q_{P}(y) := \EXP_{X \sim P}\ \big(\qyx{y}{X}\big).
\]
On its support 
\[\Y_{P} := \lp\{y \in \Y:\ q_{P}(y)>0\rp\},\]
the "proto" conditional mean $f_P(y)$ is defined as the function
\beq\label{eq:DefProtoMMSE}
y\in\Y_{P}\mapsto f_{P}(y)
:=  
\frac{\EXP_{X \sim P}\ \big(X \qyx{y}{X}\big)}{ \EXP_{X \sim P}\ \big( \qyx{y}{X} \big)}.
\eeq
\ED

\subsection{Scalar denoising}
 
In the scalar case, we fully characterize proto conditional distributions $\qyx{y}{x}$ such that the conditional mean estimator $f_P$ defined by \eq{eq:DefProtoMMSE} is a proximity operator.
\BL\lab{ScalMMSE}
Consider $\X,\Y \subset \RR$ and  $q: \X \times \Y \to \RR_+$. For $P$ a probability distribution on $\X$ such that \eq{eq:admrv} holds, let 
$f_P$, $\Y_P$ be defined as in \defi{def:ProtoMMSE}. The following properties are equivalent:
\bena
\item \lab{scalmmsea} For any $P$ satisfying~\eq{eq:admrv}, $f_P$ \rev{implements} a proximity operator;
\item \lab{scalmmseb} For any $x,x' \in \X, y,y' \in \Y$ with $x'>x$ and $y'>y$,
\(
\qyx{y'}{x'}\qyx{y}{x}-\qyx{y'}{x}\qyx{y}{x'} \geq 0.
\)
\item[]\hspace{-0.6cm}If $\qyx{y}{x}>0$ on $\X \times \Y$ then \eq{scalmmsea}-\eq{scalmmseb} are further equivalent to:
\item\lab{scala} For any $y,y' \in \Y$ with $y'>y$, the function $x \mapsto \log \qyx{y'}{x}-\log \qyx{y}{x}$ is non-decreasing;
\item\lab{scalb} For any $x,x' \in \Y$ with $x'>x$, the function $y \mapsto \log \qyx{y}{x'}-\log \qyx{y}{x}$ is non-decreasing.
\item[] \hspace{-0.6cm}If $\tfrac{\d}{\d y} \log \qyx{y}{x}$ exists on $\X \times \Y$, then \eq{scalmmsea}-\eq{scalmmseb}  are further equivalent to:
\item\lab{scald} For any $y \in \Y$, the function $x \mapsto \tfrac{\d}{\d y} \log \qyx{y}{x}$ is non-decreasing.
\item[] \hspace{-0.6cm}If $\tfrac{\d}{\d x} \log \qyx{y}{x}$ exists on $\X \times \Y$, then \eq{scalmmsea}-\eq{scalmmseb} are further equivalent to:
\item\lab{scale} For any $x \in \X$ the function $y \mapsto \tfrac{\d}{\d x} \log \qyx{y}{x}$ is non-decreasing.
\een
\EL
The proof is postponed to Annex~\ref{ann:ScalMMSE}. Two applications are scalar Poisson denoising (\prop{scalpoisson}) and denoising in the presence of additive noise (\prop{scallogconc}).

\proof[Proof of \prop{scalpoisson}]
Consider $\X :=  \RR_+^*$, $\Y := \RR_+$, and the proto-conditional distribution $\qyx{y}{x} := \frac{x^{y}}{\Gamma(y+1)} e^{-x}$ which is defined and strictly positive on $\X \x \Y$. For any $y \in \Y$ we have \(\sup_{x \in \X} (1+|x|)\ q(y|x) < \infty,\) hence property~\eq{eq:admrv} holds for any distribution $P$.
\defi{def:ProtoMMSE} yields $q_P(y)>0$ on $\Y$ hence $\Y_P = \Y$. Setting $P = p_X$, as $p_{Y|X}(Y=y|x) = \qyx{y}{x}$ for any $y \in \Y' := \NN$ we get $\EXP\ (X | Y=y) = f_{P}(y)$ with $f_{P}$ given by~\eq{eq:DefProtoMMSE}. 
On $\X \times \Y$ we have $\log \qyx{y}{x} = y \log x - \log \Gamma(y+1) -x$, hence $y \mapsto \frac{\d}{\d x} \log \qyx{y}{x} = \frac{y}{x}-1$ is non-decreasing on $\X$. The fact that $f_P$ \rev{implements} a proximity operator, i.e., the existence of $\tilde{\ph}_X$, follows by \lem{ScalMMSE}\eq{scalmmsea}$\Leftrightarrow$\eq{scald}. 
\endproof

\begin{proof}[Proof of \prop{scallogconc}]
Consider $\X = \Y = \RR$, $\qyx{y}{x} := p_{Y|X}(y|x) = \exp(-F(y-x))$, and observe that $\log \qyx{y}{x} = -F(y-x)$. For $P:=p_X$ a probability distribution satisfying \eq{eq:admrv}, reasoning as in the proof of \prop{scalpoisson} we have $q_P(y)>0$ on $\Y$ hence $\Y_P=\Y$ and  $\EXP(X|Y=y)=f_P(y)$.

Consider first the case where $F$ is $C^1$. Then, $F$ is convex iff $u \mapsto F'(u)$ is non-increasing, which is equivalent to $y \mapsto \frac{\d}{\d x} \log \qyx{y}{x} = F'(y-x)$ being non-decreasing. By \lem{ScalMMSE}\eq{scale}$\Leftrightarrow$\eq{scalmmsea} this is equivalent to the fact that $f_P$ \rev{implements} a proximity operator for any $P$ satisfying \eq{eq:admrv}.
To conclude, we exploit a simple lemma proved in Appendix~\ref{pfcvxbnd}.
\BL \lab{cvxbnd} 
If $G: \RR \to \RR$ is convex and satisfies $\int_\RR e^{-G(x)} dx < \infty$ then 
\(
\sup_{x \in \RR} (1+|x|) e^{-G(x)} < \infty.
\)
\EL
Hence:
\bit 
\item if $F$ is not convex then there exists $P$ satisfying \eq{eq:admrv} such that $f_P$ \rev{does {\em not} implement} a proximity operator, hence \prop{scallogconc}\eq{scallcb} cannot hold;
\item if $F$ is convex then for any $y$ the function $G: x \mapsto F(y-x)$ is convex and $\int_\RR e^{-G(x)}dx < \infty$ hence, by \lem{cvxbnd}, the function $x \mapsto (1+|x|) \qyx{y}{x}$ is bounded. As a result, \eq{eq:admrv} holds for any distribution $P$. As just shown, $f_P$ \rev{implements} a proximity operator {\em for any $P$}.
\eit
This establishes the equivalence \prop{scallogconc}\eq{scallca}$\Leftrightarrow$\eq{scallcb} when $F$ is $C^1$.

To extend the result when $F$ is only $C^0$, we reason similarly using \lem{ScalMMSE}\eq{scalb}$\Leftrightarrow$\eq{scalmmsea}.  A bit more work is needed to show that \lem{ScalMMSE}\eq{scalb} holds iff $F$ is convex, as we now establish. 

With the change of variable $u = y-x'$, $h = x'-x>0$, \lem{ScalMMSE}\eq{scalb} is equivalent to: 
\beq\lab{tmpscallc}
\all h>0,\qu u \mapsto F(u+h)-F(u)\qu \text{is non-decreasing}.
\eeq
When \eq{tmpscallc} holds we have for any $u_1<u_2$ (using $h := (u_2-u_1)/2$ in \eq{tmpscallc} )
\[
F(\tfrac{u_1+u_2}{2})-F(u_1) = F(u_1+h)-F(u_1) \leq F(\tfrac{u_1+u_2}{2}+h)-F(\tfrac{u_1+u_2}{2}) = F(u_2)-F(\tfrac{u_1+u_2}{2})
\]
hence $F(\tfrac{u_1+u_2}{2}) \leq \tfrac{F(u_1)+F(u_2)}{2}$. As $F$ is $C^0$, this is well known to imply that $F$ is convex. 

Vice-versa, when $F$ is convex, given $u$ and $h,h'>0$ we wish to prove that $F(u_4)-F(u_3) \geq F(u_2)-F(u_1)$ where $u_1:=u$, $u_2:=u+h$,$u_3:=u+h'$, $u_4:=u+h+h'$. Two cases are possible: either $u_1 < u_2 < u_3 < u_4$ or $u_1 < u_3 \leq  u_2 < u_4$. We treat the latter, the first one can be handled similarly. Since $F$ is convex and $u_3 \leq u_2$, there exists $a \in \d F(u_3)$ and $b \in \d F(u_2)$ with $a \leq b$, and we get $a(u_3-u_1) = ah' \leq b h' = b(u_4-u_2)$. As a result
\begin{eqnarray*}
F(u_4)-F(u_3) &=& F(u_4)-F(u_2)+F(u_2)-F(u_3)
\geq b(u_4-u_2) + F(u_2)-F(u_3) \\
&\geq& a(u_3-u_1) + F(u_2)-F(u_3)\\
&\geq& F(u_3)-F(u_1) + F(u_2)-F(u_3) = F(u_2)-F(u_1).
\end{eqnarray*}
\end{proof}

\BE[Completely separable model]\lab{compsep}
Consider a completely separable model: the entries $X_i \in \RR$, $1 \leq i \leq n$ of the random variable $X \in \RR^n$ are drawn independently with prior distributions $P_i$, i.e., $p_X(x) = \prod_{i=1}^n p_{X_i}(x_i)$; the conditional distribution of $Y \in \RR^n$ given $x \in \RR^n$ corresponds to independent noise on each coordinate, $p_{Y|X}(y|x) = \prod_{i=1}^n p_{Y_i|X_i}(y_i|x_i)$. This implies that the conditional expectation can be computed coordinate by coordinate. Assuming further that $p_{Y_i|X_i}(y_i|x_i) \propto q_i(y_i|x_i)$ where $q_i$ satisfies Lemma~\ref{ScalMMSE}-\eq{scalmmseb}, we obtain for $1 \leq i \leq n$
\[
\EXP(X|Y=y)_i = \EXP(X_i|Y_i=y_i) \rev{~\in~} \arg\min_{x_i \in \RR} \lp\{\tfrac{1}{2}(y_i-x_i)^2+\ph_{P_i}(x_i)\rp\},\qu \all y \in \RR.
\]
and as a result 
\[
\EXP(X|Y=y) \rev{~\in~} \arg\min_{x \in \RR^n} \lp\{\tfrac{1}{2}\|y-x\|^2+\sum_{i=1}^n \ph_{P_i}(x)\rp\},\qu \all y \in \RR^n.
\]
Hence, the conditional mean \rev{implements} the proximity operator of $\ph_P(x) := \sum_{i=1}^n \ph_{P_i}(x_i)$.
\EE

\subsection{Multivariate denoising}

In dimension $1 \leq n \rev{~<~} \infty$, we have the following result.
\BL \lab{VectMMSE}
Consider $\X,\Y \subset \H$ and  $q: \X \times \Y \to \RR_+$. For $P$ a probability distribution on $\X$ such that \eq{eq:admrv} holds, let 
$f_P$, $\Y_P$ be defined as in \defi{def:ProtoMMSE}.

Assume \rev{that for any $P$ satisfying \eq{eq:admrv}, $f_P$ implements a proximity operator}. Then, for any $x,x' \in \X$, $y,y' \in \Y$:
\bena
\item \lab{qneca} if $\nabla_{x}\log \qyx{y}{x}$ and $\nabla_{x}\log \qyx{y'}{x}$ exist, there is a scalar $c = c(x,y,y') \geq 0$ such that
\beq\lab{qnecaeq}
\nabla_{x} \log \qyx{y'}{x} - \nabla_{x} \log \qyx{y}{x} = c  \ (y'-y).
\eeq
\item \lab{qnecabis} if $\nabla_{y}\log \qyx{y}{x}$ and $\nabla_{y}\log \qyx{y}{x'}$ exist, there is a scalar $c = c(x,x',y) \geq 0$ such that
\beq\lab{qnecabiseq}
\nabla_{y} \log \qyx{y}{x'} - \nabla_{y} \log \qyx{y}{x} = c ' \ (x'-x).
\eeq
\een
\EL
The proof is postponed to Annex~\ref{ann:VectMMSE}. 
{\remark In \eq{qneca}-\eq{qnecabis} the gradients are only assumed to exist at particular points $x,x',y,y'$, leading to the necessary conditions \eq{qnecaeq}-\eq{qnecabiseq} at these points. }

A first consequence is that MMSE estimation with additive Laplacian noise (resp. with Poisson noise) in dimension $n \geq 2$ behaves differently from dimension $n = 1$ (cf  \prop{scalpoisson}, \prop{scallogconc}).

\begin{example}[multivariate additive Laplacian noise]\lab{exlapv}
Consider a multivariate Laplacian noise model: given $x \in \X = \RR^n$, the conditional probability of the random vector $Y$ on $\Y = \RR^{n}$ is defined by
$p(y|x) \propto \qyx{y}{x} := e^{-\|x-y\|_{1}}$. As $\log \qyx{y}{x} = -\|x-y\|_{1}$, given $y = (y_{i})_{i=1^{n}}$ and $x = (x_{i})_{i=1}^{n}$ such that $x_{i} \neq y_{i}$, $\all 1 \leq i \leq n$, $\log q(\cdot,y)$ is differentiable at $x$ and $\nabla_{x} \log \qyx{y}{x} = -\sign(x-y)$. Hence, for $x \in \X$, $y,y' \in \Y$ such that $x_{i} \notin \{y_{i},y'_{i}\}$ for $1 \leq i \leq n$ we have
\beq \lab{exlapveq}
\nabla_{x} \log \qyx{y'}{x}-\nabla_{x} \log \qyx{y}{x} = - \big(\sign(x-y')-\sign(x-y)\big).
\eeq
\bit
\item The scalar case $n=1$ is covered by \prop{scallogconc}  since $q$ is log-concave. For any distribution $P$ on the scalar random variable $X \in \RR$, the MMSE estimator $f_P$ \rev{implements} a proximity operator.
\item For $n \geq 2$, this is no longer the case. Consider $x = 0$ and $y \in (\RR_+^*)^n$ such that $y_{1} \neq y_{2}$, and $y'= -y$. The vector $-(\sign(x-y')-\sign(x-y)) = 2 \cdot \mathbf{1}_{n}$ is not proportional to the vector $(y'-y) = 2y'$ which first two entries are distinct. Hence \eq{exlapveq} is incompatible with condition~\eq{qnecaeq} and by \lem{VectMMSE}, there exists a prior distribution $P$ such that $f_P$ \rev{does {\em not} implement} a proximity operator.
\eit
\end{example}
{\remark \lab{sepvsnonsep} As the noise distribution is separable (multivariate Laplacian noise corresponds to i.i.d. scalar noise on each coordinate), by \exa{compsep} the MMSE estimator \rev{in fact implements} a proximity operator as soon as the prior $P$ is also separable (i.e., if $X$ has independent entries $X_i$). However, for $n \geq 2$, we have just shown that there exists a prior $P$ (non-separable) such that $f_P$ \rev{does not implement} a proximity operator.
}

\begin{example}[multivariate Poisson denoising]\lab{poisson}
Consider a multivariate Poisson noise model: given $x \in \X := (\RR_+^*)^{n}$, the conditional probability of the random vector of integers $Y$ on $\Y := \NN^{n}$ is defined by
$p(y|x) = \qyx{y}{x} := \prod_{i=1}^{n} \left(\frac{x_{i}^{y_{i}} e^{-x_{i}}}{\Gamma(y_{i}+1)}\right)$. We observe that
\[
\log \qyx{y}{x} = \sum_{i=1}^{n} y_{i} \ln x_{i} -x_{i} - \log \Gamma(y_{i}+1)
\]
Given $y \in \Y$, the function $x \mapsto \log \qyx{y}{x}$ is differentiable on $\X$ with $\nabla_{x} \log \qyx{y}{x} = (y_{i}/x_{i}-1)_{i=1}^{n}$ hence for $x \in \X$, $y,y' \in \Y$ we have
\[
\nabla_{x} \log \qyx{y'}{x}-\nabla_{x} \log \qyx{y}{x} = ((y'_{i}-y_{i})/x_{i})_{i=1}^{n}.
\]
\bit 
\item The case $n=1$ is covered by \prop{scalpoisson}: $f_P$ \rev{implements} a proximity operator for any prior $P$.
\item For $n \geq 2$ this is again no longer the case. Consider, e.g., $x \in \X$ with $x_{1} \neq x_{2}$ and $y,y' \in \Y$ such that $y_i \neq y'_i$, $i=1,2$. The vector \rev{with entries} $(y'_i-y_i)/x_i$ cannot be proportional to \rev{the vector} $y'-y$ hence \eq{qnecaeq} cannot hold. By \lem{VectMMSE}, there is a prior $P$ such that $f_P$ \rev{does {\em not} implement} a proximity operator. 
\eit
\end{example}
\rem{sepvsnonsep} again applies: for a separable prior $P$, $f_P$ \rev{implements} a proximity operator, yet there exists a (non-separable) prior $P$ such that $f_P$ \rev{does not implement} a proximity operator.

For smooth enough proto-conditional distributions we have the following corollary of \lem{VectMMSE}.

\BL \lab{VectMMSENec}
Let $\H$ be of dimension $2 \leq n \leq \infty$ and consider open sets $\X,\Y \subset \H$ where $\X$ is connected. Let $q: \X \times \Y \to \RR_+^*$ be such that $x \mapsto \nabla_x \log \qyx{y}{x}$ is $C^0$ for any $y \in \Y$, and $y \mapsto \nabla_y \log \qyx{y}{x}$ is $C^0$ for any $x \in \X$. For $P$ 
a probability distribution on $\X$ satisfying \eq{eq:admrv}, let $f_P$, $\Y_P$ be defined as in \defi{def:ProtoMMSE}. 

If $f_P$ \rev{implements} a proximity operator for any such $P$, then
there exists $c \geq 0$, $a \in C^1(\X)$ $b \in C^1(\Y)$ such that
\beq\lab{expfam}
\qyx{y}{x} = \exp\lp(-a(x)-b(y) + c \< x,y\>\rp),\qu \all x,y \in \X \times \Y.
\eeq
\EL
The proof is postponed to Annex~\ref{PfVectMMSENec}.
A converse result also holds.

\BL \lab{VectExpFam} 
Consider $\X,\Y \subset \H$ where $\Y$ is open and convex. Assume that $q: \X \times Y \to \RR_+^*$ has the expression \eq{expfam} where $c \geq 0$, $a: \X \to \RR$, $b: \Y \to \RR$, $b$ is $C^1(\Y)$, and $b'$ is locally bounded\footnote{If $\H$ \rev{is finite dimensional then} local boundedness is automatic by compactness arguments; the assumption is useful in infinite dimension.}. 

Let $P$ be a probability distribution satisfying: for any $y \in \Y$, there is $r = r(y)>0$ such that 
\beq \lab{admPExpFam}
\EXP_{X \sim P}\ \{(1+\|X\|)^2 \cdot \exp\lp(-a(X)+cr\|X\|+c\<X,y\>\rp)\} < \infty.
\eeq
Then \eq{eq:admrv} holds and, using the notations of \defi{def:ProtoMMSE}, we have 
\bit
\item the proto-conditional mean $f_{P}$ is differentiable on $\Y_P=\Y$;
\item its differential is symmetric positive semi-definite\footnote{
\rev{see footnote ~\ref{foot1} page~\pageref{foot1}.}
}, i.e., 
\(
Df_{P}(y) \succeq 0, \all y \in \Y;
\)
\item the proto-conditional mean $f_P$ \rev{implements the proximity operator of a penalty} $\ph_P$.
\eit
Moreover if $c>0$ and if the support of the distribution $P$ is not included in any hyperplane (i.e. if $P(\<X,u\>=d) < 1$ for any nonzero $u \in \H$ and any $d \in \RR$), then $\ph_P$ can be chosen such that:
\bit\item $Df_P(y) \succ 0$ for any $y$;
\item the proto-conditional mean $f_P$ is injective;
\item for any $y$, $f_P(y)$ is the unique stationary point (and global minimizer) of $x \mapsto \tfrac{1}{2}\|y-x\|^2+\ph_P(x)$. \rev{In particular, $f_{P}$ {\em is} the proximity operator of $\ph_{P}$.}
\eit
\EL
The proof is postponed to Appendix~\ref{pfVectExpFam}.

{\remark The case $c=0$ corresponds to independent random variables $X$ and $Y$ governed by $p_{Y|X}(y|x) = \tfrac{\qyx{y}{x}}{\int \qyx{y'}{x} dy'} = \tfrac{\exp(-a(x)-b(y))}{\int \exp(-a(x)-b(y'))dy'} = \tfrac{\exp(-b(y))}{\int \exp(-b(y'))dy'}$.}

A consequence of \lem{VectExpFam} is that we recover the results of \cite{GRIBONVAL:2010:INRIA-00486840:1}. We also cover a variant of multivariate Poisson denoising that leads again to a proximity operator even for $n \geq 2$ as expressed in \prop{logpoisson}.

\subsubsection{Proof of \prop{gauss}}
Consider white Gaussian noise, i.e. $p(y|x) \propto \qyx{y}{x}$ with $\qyx{y}{x} := \exp\lp(-\tfrac{c}{2} \|x-y\|^2\rp)$; $\X = \Y = \H = \RR^{n}$. Observe that $\log \qyx{y}{x} = -\tfrac{c}{2} \|x-y\|^2 = c\<x,y\>-a(x)-b(y)$ with $a(x):=\tfrac{c}{2}\|x\|^2$, $b(y) := \tfrac{c}{2}\|y\|^2$.
Since $\sup_{x \in \X} (1+\|x\|)^2 \exp\lp(-a(x)+(r+c\|y\|)\|x\|\rp) \leq \sup_{u \in \RR_+} (1+u)^2 e^{-cu^2/2+(r+c\|y\|)u} < \infty$, any probability distribution $P$ satisfies \eq{admPExpFam} hence
we can apply \lem{VectExpFam}.

\subsubsection{Proof of \prop{logpoisson}}

Consider $\X := (\RR_+^*)^n$, $\Y' = \NN^n$, $\Y := (\RR_+)^n$. For $(x,y) \in \X \times \Y'$ we have $p(y|x) = \qyx{y}{x}$ where
\[
\log \qyx{y}{x} = \sum_{i=1}^{n} \lp(y_i \log x_{i} -x_{i} - \log \Gamma(y_{i}+1)\rp)
\]
is defined on $\X \times \Y$.
Denoting $z := (\log x_i)_{i=1}^n \in \Z := \RR^n$ we have $\wt{p}(y|z) := \wt{q}(y|z)$ where $\log \wt{q}(y|z) := \<z,y\>-a(z)-b(y)$ with \[
a(z):=\sum_{i=1}^n e^{z_i}
\]
and $b(y):= \sum_{i=1}^n \log \Gamma(y_i+1)$.  

For any $y \in \Y$, $r>0$, using arguments similar to those in the proof of \lem{cvxbnd}, we get that
\[
\sup_{z \in \Z} (1+\|z\|)^2 \exp\lp(-a(z)+(r+\|y\|)\|z\|\rp)
< \infty
\]
we get that for any distribution $p_X$ on $X$, denoting $P$ the resulting distribution on the random variable $Z = \log X \in \Z$, $P$ necessarily satisfies \eq{admPExpFam} for any $y \in \Y$ and $r>0$. 

As $b$ is $C^1(\Y)$, both $b$ and $b'$ are locally bounded hence, we can again apply \lem{VectExpFam} to get $ f_P(y) = \EXP_{Z \sim P} (Z | Y=y) = \EXP (\log X|Y=y)$ \rev{implements} a proximity operator. 

\appendix
\section{Proofs} 

\subsection{Proof of \lem{ScalMMSE}} \lab{ann:ScalMMSE}

\eq{scalmmsea}$\Rightarrow$\eq{scalmmseb}. First, we establish a necessary condition that any function $\qyx{x}{y}$ must satisfy to ensure that \eq{eq:DefProtoMMSE} \rev{implements} a proximity operator for any prior probability distribution $P$ on the random variable $X$. This condition will be re-used for the proof of \lem{VectMMSE}.

\BL \lab{lem:ProxMMSENec}
Consider $\X,\Y \subset \H$ and $q : \X \times \Y \to \RR_+, (x,y) \mapsto \qyx{y}{x}$. Assume that {\em for any probability distribution $P$ such that~\eq{eq:admrv} holds}, $f_P$ \rev{implements} the proximity operator of some penalty $\ph_P$ (i.e. $f_P(y) \in \prox_{\ph_P}(y)$, $\all y \in \Y_P$), with $\Y_P$ and $f_P$ as in \defi{def:ProtoMMSE}. Then the function $q$ satisfies
\beq\label{eq:NSCProtoMMSE}
\all x,x' \in \X, y,y' \in \Y,\qu\qu
\<x'-x,y'-y\> \left(\qyx{y'}{x'}\qyx{y}{x}-\qyx{y'}{x}\qyx{y}{x'}\right) \geq 0.
\eeq
\EL 
Even though \eq{eq:NSCProtoMMSE} is not intuitive, its main interest is that it depends only on the ``proto'' conditional distribution $\qyx{y}{x}$ {\em but not on the prior distribution $P$ on $X$}. It \rev{necessarily} holds when the "proto" conditional distribution $\qyx{y}{x}$ \rev{is such} that {\em for any prior probability distribution $P$ on $X$}, the "proto" conditional expectation $f_P$ \rev{implements} a proximity operator. 

\proof
Given $x,x' \in \X$, $y,y'\in \Y$, consider the probability distribution $P := \tfrac{1}{2}( \delta_{x}+\delta_{x'})$. It is straightforward that $P$ satisfies 
\eq{eq:admrv}.

Now, denoting $\lambda := (\qyx{y}{x}+\qyx{y}{x'})(\qyx{y'}{x}+\qyx{y'}{x'})$, we distinguish two cases depending whether $\lambda = 0$ or not.

First, if $\lambda = 0$ then one must have $\qyx{y}{x}+\qyx{y}{x'}=0$, or $\qyx{y'}{x}+\qyx{y'}{x'}=0$, or both. 
Without loss of generality we treat the case $\qyx{y}{x}+\qyx{y}{x'}=0$ (the other cases are treated similarly). Since $q$ is non-negative, this implies $\qyx{y}{x}=\qyx{y}{x'}=0$, hence $\qyx{y'}{x'}\qyx{y}{x}-\qyx{y'}{x}\qyx{y}{x'}= 0$ and \eq{eq:NSCProtoMMSE} trivially holds.

Now, consider the case $\lambda > 0$, i.e. $\qyx{y}{x}+\qyx{y}{x'}>0$ and $\qyx{y'}{x}+\qyx{y'}{x'}>0$. By the definition of $q_P$ and $\Y_P$ in \defi{def:ProtoMMSE}, this implies 
$y,y' \in \Y_{P}$ and thus  $q_P(y) >0$ and $q_P(y')>0$. Hence,  by the assumption of \lem{lem:ProxMMSENec}, $f_P$ \rev{implements} a proximity operator, therefore by \thm{ProCon} there exists a convex lsc function $\psi_P$ such that $f_{P}(y) \in \d\psi_P(y)$ and $f_{P}(y') \in \d\psi_P(y')$. From the definition of a subdifferential, this implies that $\psi_P(y')-\psi_P(y) \geq \<f_{P}(y),y'-y\>$ and that $\psi_P(y)-\psi_P(y') \geq \<f_{P}(y'),y-y'\>$. Adding both inequalities yields\footnote{or equivalently $\<f_{P}(y')-f_{P}(y),y'-y\> \geq 0$. Since this holds for any $y,y' \in \Y$, $f_{P}$ is a {\em monotone} operator (see, e.g., \cite{Moreau65}).}
\(
0 \geq \<f_{P}(y)-f_{P}(y'),y'-y\>.
\)
Since by \eq{eq:DefProtoMMSE}
\[
f_{P}(y) = \frac{x \qyx{y}{x}+x' \qyx{y}{x'}}{\qyx{y}{x}+\qyx{y}{x'}}\qu\mb{and}\qu f_{P}(y') = \frac{x \qyx{y'}{x}+x' \qyx{y'}{x'}}{\qyx{y'}{x}+\qyx{y'}{x'}}
\]
we obtain with straightforward computations that
\[
0 \leq \lambda \<f_{P}(y')-f_{P}(y),y'-y\> =  \<x'-x,y'-y\> \left(\qyx{y'}{x'}\qyx{y}{x}-\qyx{y'}{x}\qyx{y}{x'}\right).
\]
\endproof

In the scalar case \eq{eq:NSCProtoMMSE} is equivalent to \lem{ScalMMSE}\eq{scalmmseb}, hence \lem{lem:ProxMMSENec} establishes \eq{scalmmsea}$\Rightarrow$\eq{scalmmseb}. 

\eq{scalmmseb}$\Rightarrow$\eq{scalmmsea}. Since we consider the scalar case, the assumption that \eq{scalmmseb} holds implies \eq{eq:NSCProtoMMSE}. Consider $y,y' \in \Y_{P}$ and $X,X'$ two independent random variables with distribution $P$. Write
\beqn
q_P(y)q_{P}(y')\ \big(f_P(y')-f_P(y)\big)  
&\stackrel{\eq{eq:DefProtoMMSE}}{=}&
 q_{P}(y)\ \EXP_{X'} \lp(X'\,  \qyx{y'}{X'}\rp)\ 
 -q_{P}(y')\ \EXP_{X'} \lp(X'\, \qyx{y}{X'}\rp)
\notag\\
  &=&
\EXP_{X} \big(\qyx{y}{X}\big)\ \EXP_{X'}  \lp(X'\,\qyx{y'}{X'}\rp)\notag\\
&& -\EXP_{X}\lp(\qyx{y'}{X}\rp)\ \EXP_{X'}  \lp(X'\,\qyx{y}{X'}\rp)
  \notag\\
&=&
\EXP_{X,X'}  \Big(X' \left(\qyx{y}{X}\qyx{y'}{X'}-\qyx{y'}{X}\qyx{y}{X'}\right)\Big)
  \lab{tmp1}\\
  &\stackrel{X\leftrightarrow X'}{=}&
\EXP_{X',X}  \Big(X \left(\qyx{y}{X'}\qyx{y'}{X}-\qyx{y'}{X'}\qyx{y}{X}\right)\Big)
\lab{tmp2}\\
  & \stackrel{(\eq{tmp1}+\eq{tmp2})/2}{=}&
   \EXP_{X',X}  \Big(\tfrac{X'-X}{2} \lp(\qyx{y'}{X'}\qyx{y}{X}-\qyx{y'}{X}\qyx{y}{X'}\rp)\Big) \notag
\eeqn
If follows that
\beqnn
(y'-y)\big(f_{P}(y')-f_{P}(y)\big) &=& \tfrac{1}{2q_{P}(y)q_{P}(y')} \cdot \EXP_{X',X}\ \left\{\tfrac{(X'-X)(y'-y)}{2} \left(\qyx{y'}{X'}\qyx{y}{X}-\qyx{y'}{X}\qyx{y}{X'}\right) \right\}\\
&\stackrel{\eq{eq:NSCProtoMMSE}}{\geq}&  0\notag,
\eeqnn
hence $f_{P}$ is non-decreasing on 
the set $\Y_{P}$. By \cor{coscal}, $f_P$ \rev{implements} a proximity operator.

\eq{scalmmseb}$\Leftrightarrow$\eq{scala}$\Leftrightarrow$\eq{scalb} when $\qyx{y}{x}>0$ on $\X \times \Y$. We sketch the proof of the equivalence of~\eqref{scalmmseb} and \eq{scala}. The equivalence between~\eqref{scalmmseb} and \eq{scalb} follows similarly. Denote $Q(x;y,y'):= \log \qyx{y'}{x}-\log \qyx{y}{x}$. Property~\eq{scalmmseb} holds if and only if for any $x,x' \in \X$, $y,y' \in \Y$ with $x'>x$, $y'>y$ we have $\qyx{y'}{x'} \qyx{y}{x} \geq \qyx{y'}{x} \qyx{y}{x'}$, which is equivalent to $\log \qyx{y'}{x'}+\log \qyx{y}{x} - \log \qyx{y'}{x}-\log \qyx{y}{x'} \geq 0$, that is to say 
\(
Q(x';y,y')-Q(x;y,y') \geq 0. 
\)
i.e., $x \mapsto Q(x;y,y')$ is non-decreasing as soon as $y'>y$. 

\eq{scala}$\Leftrightarrow$\eq{scald} and \eq{scalb}$\Leftrightarrow$\eq{scale}. This is straightforward 
using, e.g., that $\tfrac{\partial}{\partial y} \log \qyx{y}{x} = \lim_{y' \to y} \tfrac{Q(x;y,y')}{y'-y}$.

\subsection{Proof of \lem{VectMMSE}} \lab{ann:VectMMSE}

We establish \eq{qneca}.\eq{qnecabis} is obtained similarly by reversing the roles of $x$ and $y$.

Since $f_P$ is a proximity operator for any probability distribution $P$ such that \eq{eq:admrv} holds, by \lem{lem:ProxMMSENec} it follows that property \eq{eq:NSCProtoMMSE} holds. 

Consider $h \in \H$ such that $\<y'-y,h\> \dst{>} 0$. As $\log q(\cdot,y)$ is differentiable at $x$, we have $x':=x+\varepsilon h \in \X$  for $\varepsilon>0$ small enough, and $\<x'-x,y'-y\> = \varepsilon \<h,y'-y\> > 0$. By~\eq{eq:NSCProtoMMSE} we have $\qyx{y'}{x'}\qyx{y}{x} \geq \qyx{y'}{x}\qyx{y}{x'}$. Taking the logarithm yields
\[
\log \qyx{y'}{x+\varepsilon h}-\log \qyx{y'}{x} \geq \log \qyx{y}{x+\varepsilon h}-\log \qyx{y}{x}
\]
As $\log q(\cdot,y)$  and $\log q(\cdot,y')$ are both differentiable at $x$, it follows that
\[
\lp\<\nabla_{x}\log \qyx{y'}{x}, \varepsilon h\rp\> + o(\varepsilon) \geq \lp\<\nabla_{x}\log \qyx{y}{x},\varepsilon h\rp\> + o(\varepsilon)
\]
hence $\< \nabla_{x}\log \qyx{y'}{x}-\nabla_{x} \log \qyx{y}{x},h\> \geq 0$.
Since this holds for any $h$ such that $\<y'-y,h \> \dst{>} 0$, there is a scalar $c \geq 0$ (a priori dependent on $x,y,y'$)  such that \eq{qnecaeq} holds.

 \subsection{Proof of \lem{VectMMSENec}} \lab{PfVectMMSENec}

Since $\dim \H \geq 2$ and $\X$ and $\Y$ are open, the affine dimension of $\X$ (resp. of $\Y$) exceeds two. We use the following lemma.
\BL \lab{lemaff0}
Let $\Z \subset \H$ be of affine dimension at least two (i.e., there is no pair $z_1,z_2 \in \H$ such that $\Z \subset \{t z_1 + (1-t) z_2,\ t \in \RR\}$). Assume that the function $\ff: \Z \to \H$ satisfies 
\beq \lab{aff0}
\all z,z' \in \Z,\ \exists c(z,z') \in \RR,\ \ff(z)-\ff(z') = c(z,z')\ (z-z').
\eeq
Then, there exists $c \in \RR$ 
such that
\[
\all z,z' \in \Z,\ \ff(z)-\ff(z') = c\ (z-z').
\]
\EL
\proof
{\bf Step 1.} We show that if $z_1,z_2,z_3 \in \Z$ are not aligned (i.e., if they are affinely independent) then $c(z_i,z_j)=c(z_1,z_2)$ for all $1 \leq i \neq j \leq 3$.

Denote $c_{ij} := c(z_i,z_j)$. By symmetry\eq{aff0} yields $c_{ij}=c_{ji}$, $\all i \neq j$. Moreover, summing up \eq{aff0} with $z=z_i$, $z'=z_j$ over $(i,j) \in \{(1,2),(2,3),(3,1)\}$ yields
\[
0 = c_{12}(z_{1}-z_{2})+c_{23}(z_{2}-z_{3}) + c_{31}(z_{3}-z_{1}) 
   = (c_{12}-c_{31})z_{1} + (c_{23}-c_{12}) z_{2}  + (c_{31}-c_{23}) z_{3}.
\]
As the coefficients in the right hand side sum to zero, affine independence yields $c_{12}=c_{23}=c_{31}$.

{\bf Step 2.}
As the affine dimension of $\Z$ exceeds two it is not a singleton and we can choose an arbitrary pair $z_1 \neq z_2 \in \Z$. Define $c := c(z_1,z_2)$. We show that for any $x, y \in \Z$ we have $c(x,y)=c(y,x)=c$.

Define $\S := \{z_1+(1-t)z_2,\ t \in \RR\}$ the affine hull of $z_1,z_2$ and observe that $\Z \cap \S^c \neq \void$.

First, consider $x \in \Z \cap \S^c$. For $y=z_1$, as $z_1,z_2,x$ are not aligned, by Step 1 we have $c(x,z_1) = c(z_1,x) = c$. For $y \in \Z \cap \S \backslash \{z_1\}$, as $x,y,z_1$ are not aligned, by Step 1 again we get $c(x,y) = c(y,x) = c(x,z_1)=c$. This establishes the result for any $x \in \Z \cap \S^c$ and $y \in \Z \cap \S$.

Second, consider $x, y \in \Z \cap \S^c$. As just shown, we have $c(x,z_1) = c(z_1,y) = c$, hence by \eq{aff0}
\[
\ff(x)-\ff(y) = \ff(x)-\ff(z_1)+\ff(z_1)-\ff(y) = c(x,z_1)\ (x-z_1)+c(z_1,y)\ (z_1-y) = c\ (x-y).
\]

Finally consider $x,y \in \Z \cap \S$. Let $z \in \Z \cap \S^c$ be arbitrary. As $c(x,z) = c(z,y) = c$, \eq{aff0} yields
\[
\ff(x)-\ff(y) = \ff(x)-\ff(z)+\ff(z)-\ff(y) = c(x,z)\ (x-z)+c(z,y)\ (z-y) = c (x-y).
\]
\endproof

For a given $x \in \X$, consider the function $y \mapsto \ff_x(y) := \nabla_x \log \qyx{y}{x}$. By \lem{VectMMSE}\eq{qneca}, $\ff_x$ satisfies \eq{aff0} and the constants $c(z,z')$ are non-negative. By \lem{lemaff0} there  
is $c_x \in \RR_+$ 
such that 
\beq \lab{vn0}
\ff_x(y)-\ff_x(y') = c_x\ (y-y'),\ \all y,y' \in \Y.
\eeq
As a result $y \mapsto \ff_x(y)$ is differentiable with $D_y \nabla_x \log \qyx{y}{x} = c_x\ \Id$.
Similarly, given $y \in \Y$, with $\gg_y(x) := \nabla_y \log \qyx{y}{x}$, there 
is $d_y \in \RR_+$ 
such that 
\beq \lab{vn1}
\gg_y(x)-\gg_y(x') = d_y\ (x-x'),\ \all x,x' \in \X
\eeq
hence $D_x \nabla_y \log \qyx{y}{x} = d_y\ \Id$.

When $(x,y) \mapsto \log \qyx{y}{x}$ is $C^2$, by Schwarz' theorem we have $D_x \nabla_y \log \qyx{y}{x} = D_y \nabla_x \log \qyx{y}{x}$ for any $x,y \in \X \times \Y$. Thus, $c_x\ \Id = d_y\ \Id$ for any $x,y \in \X \times \Y$ hence $c_x=d_y=c\geq 0$ is independent of $x,y$ and  
\beq \lab{vn2}
\nabla_x \log \qyx{y}{x}-\nabla_x \log \qyx{y'}{x} = c (y-y'),\qu \all x \in \X, \all y,y' \in \Y.
\eeq
Let us now show that \eq{vn2} also holds with the considered weaker assumption on $q$ . For this, fix some arbitrary $x \in \X$ and $x' \in \X$ close enough so that $\{x+t(x'-x)\} \subset \X$ (remember that $\X$ is open so this is possible). Consider $y,y' \in \Y$. Denote 
\[
f(t):= \log \qyx{y}{x+t(x'-x)}-\log \qyx{y'}{x+t(x'-x)},\ t \in [0,1].
\]
As $x \mapsto \nabla_x \log \qyx{y}{x}$ is assumed to be continuous, the function $f$ is $C^1$ on $[0,1]$ and by \eq{vn0} we have
\[
f'(t) = \< \nabla_x \log \qyx{y}{x+t(x'-x)}-\nabla_x \log \qyx{y'}{x+t(x'-x)},x'-x\> 
= c_{x+t(x'-x)} \<y-y',x'-x\>.
\]
As a result
\beqnn
\log \qyx{y}{x'}-\log \qyx{y'}{x'}
-\log \qyx{y}{x}+\log \qyx{y'}{x}
&=&
f(1)-f(0) = \int_0^1 f'(t) dt\\
& =& \int_0^1 c_{x+t(x'-x)} dt\ \<y-y',x'-x\>.,
\eeqnn
By \eq{vn1}, taking the gradient of both sides with respect to $y$ yields
\[
d_y\ (x'-x) =  \int_0^1 c_{x+t(x'-x)} dt\ (x'-x).
\]
Thus, $d_y$ does not depend on $y$. Similarly, $c_x$ does not depend on $x$. This establishes \eq{vn1}.

Fix an arbitrary $y_0 \in \Y$ and denote $H(x,y) := \log \qyx{y}{x}-c\<x,y\>$ and $a(x) := -H(x,y_0)$. We obtain
\[
\nabla_x \lp(H(x,y)+a(x)\rp) = \nabla_x \lp( H(x,y)-H(x,y_0)\rp) = 0,\qu \all x \in \X, \all y \in \Y.
\]
Since $\X$ is open and connected, hence path-connected, it follows that for any $y \in \Y$ there is $b(y) \in \RR$ such that: $H(x,y)+a(x) = -b(y)$ for all $x \in \X$, i.e.,
\[
\log \qyx{y}{x} = -a(x)-b(y)+c\<x,y\>.
\]
As both $x \mapsto \nabla_y \log \qyx{y}{x}$ and $y \mapsto \nabla_x \log \qyx{y}{x}$ are $C^1$, both $a(\cdot)$ and $b(\cdot)$ are $C^1$.

\subsection{Proof of \lem{VectExpFam}} \lab{pfVectExpFam}

We use the following lemma, which proof is slightly postponed.
\BL\label{lem:JacobianMMSE}
Consider \rev{$\H$ a Hilbert space,}  two subsets $\X,\Y \subset \H$ and a "proto" conditional distribution $q: \X \times \Y \to \RR_+$. Let $\V \subset \Y$ be an open set such that  the  gradient $\nabla_{y}\qyx{y}{x}$ exists for any $x \in \X$, $y \in \V$. Define
\beq
C_{q,\V}(x) := \sup_{y \in \V}   \left\{\qyx{y}{x} + \|\nabla_{y} \qyx{y}{x}\|_{\H}\right\} \in [0,\infty] \label{eq:QBounded}
\eeq
Consider $P$ a probability distribution such that 
\beq\lab{admrvd}
\EXP_{X \sim P}\ \{ (1+ \|X\|_{\H}) \  C_{q,\V}(X)\} < \infty.
\eeq
Then \eq{eq:admrv} holds, the set $\V_{P} := \{y \in \V:\ q_{\dst{P}}(y)>0\} = \V \cap \Y_P$ is open, and
\bit
\item the "proto" marginal distribution $q_{P}(y) := \EXP_{X \sim P} \ \lp(\qyx{y}{X}\rp)$ is differentiable on $\V$;
\item the "proto" conditional mean $f_{P}$ defined by~\eq{eq:DefProtoMMSE} is differentiable on $\V_{P}$ with\footnote{For $u,v \in \H$, $v^{T}: \H \to \RR$ denotes the linear form $x \mapsto \<v,x\>$, and $u v^{T}: \H \to \H$ the linear operator $x \mapsto \<v,x\> u$.}
\begin{equation}
Df_{P}(y) = 
\tfrac{1}{2q_{P}^2(y)} \ \EXP_{X \sim P} \EXP_{X' \sim P}  \Big(
\qyx{y}{X}
\qyx{y}{X'} \cdot (X'-X)\left(
\nabla_{y}^{T} \log \qyx{y}{X'} -  \nabla_{y}^{T} \log \qyx{y}{X}
\right)\Big)\label{eq:JacobianCharacterization}
\end{equation}
\eit
where by convention $\nabla_{y} \log \qyx{y}{x} = 0$ when $\qyx{y}{x}=0$.
\EL
{\remark {The assumption \eq{admrvd} is chosen for simplicity but could be relaxed.\\}}

By~\eq{expfam} and the fact that $b \in C^1(\Y)$ we have for all $(x,y) \in \X \times \Y$
\begin{eqnarray*}
\qyx{y}{x}+ \|\nabla_y \qyx{y}{x}\| 
&=& \qyx{y}{x} \lp(1+\|\nabla_y \log \qyx{y}{x}\|\rp)\\
&=& \exp\lp(-a(x)-b(y)+c\<x,y\>\rp)  \lp(1+\|b'(y)+x\|\rp)\\
&\leq&
\exp\lp(-a(x)+c\<x,y\>\rp) \exp\lp(-b(y)\rp) \lp(1+\|b'(y)\|\rp)\lp(1+\|x\|\rp).
\end{eqnarray*}
Consider $y \in \Y$. Since $\Y$ is open, there is $r_0$ such that for $0<r<r_0$ the closed ball $\V=B(y,r)$ is a neighborhood of $y$ satisfying $\V \subset \Y$. By the local boundedness of $b'$, $b$ is locally Lipschitz hence locally bounded hence there is $r_1$ such that $C(\V) := \sup_{y' \in \V} \exp\lp(-b(y)\rp) \lp(1+\|b'(y)\|\rp) < \infty$ for any $0<r < r_1$. For $0<r<\min(r_0,r_1)$ we have
$\sup_{y' \in \V} \exp\lp(c\<x,y'-y\>\rp) = \exp\lp(cr \|x\|\rp) < \infty$, hence for any $x \in \X$
\begin{eqnarray*}
C_{q,\V}(x) &:=& \sup_{y' \in \V} \lp\{\qyx{y'}{x}+ \|\nabla_y \qyx{y'}{x}\|\rp\}\\
&\leq & C(\V) \sup_{y' \in \V} \exp\lp(-a(x)+c\<x,y'-y\>+c\<x,y\>\rp)  \lp(1+\|x\|\rp)\\
& \leq & C(\V) \exp\lp(-a(x)+ cr \|x\|+c\<x,y\>\rp) \lp(1+\|x\|\rp).
\end{eqnarray*}
By assumption~\eq{admPExpFam} we have $\EXP_{X\sim P} \lp\{(1+\|X\|) C_{q,\V}(X)\rp\} < \infty$ when $r>0$ is small enough, i.e. \eq{admrvd} holds. Moreover, since $\qyx{y}{x}>0$ for any $x \in \X$ we have $q_P(y)>0$ i.e. $y \in \Y_P$, showing that $\Y=\Y_P$. By \lem{lem:JacobianMMSE}, $f_{P}$ is differentiable at $y$ with differential given by \eq{eq:JacobianCharacterization}. 
Finally, by \eq{expfam} we have
\[
(x'-x) \left(\nabla_{y}^{T}\log \qyx{y}{x'}-\nabla_{y}^{T}\log \qyx{y}{x}\right)
= c(x'-x)(x'-x)^{T} \succeq 0,\qu \all x,x' \in \X
\]
hence $Df_{P}(y)$ is the expectation of a symmetric positive semi-definite operator. As a result, $Df_P(y)$ is symmetric positive semi-definite. As this holds for any $y \in \Y$, and since $\Y$ is open and convex, by \thm{jac1}, $f_P$ \rev{implements} a proximity operator, i.e., there exists a function $\ph_P$ such that \rev{$f_P(y) \in \prox_{\ph_P}(y)$ for all $y \in \Y$}. 

Now, assume that $c>0$ and consider $u \in \H$ such that $\<u,Df_P(y) u\> =0$. As
\[
\<u,Df_P(y) u\> = \tfrac{c}{q_P^2(y)} \EXP_{X \sim P} \EXP_{X' \sim P} q(y|X)q(y|X') \<u,(X'-X)\>^2.
\]
and $q(y|x)q(y|x') \<u,x'-x\>^2 \geq 0$ for any $x,x'$, by Markov's inequality we obtain that 
\[
q(y|X)q(y|X') \<u,X'-X\>^2 =0
\]
almost surely on the draw of $(X,X')$. As $q(y|x)>0$ for any $x$, this implies that $\<u,X'-X\>=0$ almost surely, hence there exists $d \in \RR$ such that $\<u,X\>=d$ almost surely.
Since we assume that $P(\<u,X\>=d) < 1$ for any nonzero $u \in \H$, it follows that $u = 0$. This shows that $Df_P(y) \succ 0$. We conclude using \cor{uniqueC1}.

\subsection{Proof of \lem{lem:JacobianMMSE}} 

By~\eqref{eq:QBounded}, $0 \leq \qyx{y}{x} \leq C(x)$ for any $y \in \V,x \in \X$, hence the numerator in~\eqref{eq:DefProtoMMSE}, $n(y):= \EXP_{X \sim P}\ \{X\qyx{y}{X}\}$, is well-defined. 
As, $\|\nabla_{y}\qyx{y}{x}\|_{\H} \leq C(x)$, $\EXP_{X \sim P}\ \{X \nabla_{y} \qyx{y}{X}\} <\infty$ is similarly well-defined. Denote 
$\Delta (x,y,h):= \qyx{y+h}{x}-\qyx{y}{x}-\<\nabla_{y}\qyx{y}{x},h\>$. As $|\Delta (x,y,h)| \leq 2 C(x) \ \|h\|$ for any $x \in \X,y \in \V$ and $h$ with $\|h\|_{\H}$ small enough, and as $\lim_{\|h\| \to 0} x \Delta(x,y,h) = 0$ for any $x \in \X,y \in \V$, by the dominated convergence theorem it follows that for any $y \in \V$
\[
\lim_{\|h\|_{\H} \to 0} n(y+h)-n(y)-\EXP_{X \sim P}\ \{X \<\nabla_{y}\qyx{y}{X},h\>\} = 0
\]
showing that $n$ is differentiable on $\V$ with differential $D n(y) = \EXP_{X \sim P}\ \{X\nabla_{y}^{T}\qyx{y}{X}\}$. 
Similar arguments exploiting~\eqref{eq:QBounded} show that $q_{P}$ is differentiable on $\V$ with differential $Dq_{P}(y) = \nabla^{T}q_{P}(y)$ where $\nabla q_{P}(y) = \EXP_{X \sim P}\ \{\nabla_{y}\qyx{y}{X}\}$. 
In particular, $q_{P}$ is continuous on $\V$, hence $\V_{P} = q_{P}^{-1}( (0,\infty))$ is open.  

For $y \in \V_{P}$, the denominator in~\eqref{eq:DefProtoMMSE} is $q_{P}(y)>0$. By standard calculus $f_{P}$ is differentiable at $y$ and
\begin{eqnarray*}
Df_{P}(y)
&=&
\frac{q_P(y) Dn(y)  -n(y) Dq_P(y)}{q_P^2(y)} \\
&=& 
\frac{
\EXP_{X}\ \left\{X \nabla_{y}^{T} \qyx{y}{X}\right\}
\cdot 
\EXP_{X' }\ \left\{\qyx{y}{X'}\right\}
-
\EXP_{X}\ \left\{X \qyx{y}{X}\right\}
\cdot 
\EXP_{X'}\ \left\{\nabla_{y}^{T} \qyx{y}{X'}\right\}
}{q_{P}^2(y)}\\
q_{P}^2(y) Df_{P}(y)
&=& 
\EXP_{X,X'}\ 
\left\{
X
\left(
\qyx{y}{X'} \nabla_{y}^{T} \qyx{y}{X} -  \qyx{y}{X} \nabla_{y}^{T} \qyx{y}{X'}
\right)\right\}.
\end{eqnarray*}
Now,  we distinguish two cases:
\bit
\item for $x,x'$ such that $\qyx{y}{x}\qyx{y}{x'}>0$ using that $\nabla_{y} \log q_{P} = (\nabla_{y} q_{P})/q_{P}$ where $q_{P} >0$ we write
\begin{equation}\label{eq:Simplified}
\qyx{y}{x'} \nabla_{y}^{T} \qyx{y}{x} -  \qyx{y}{x} \nabla_{y}^{T} \qyx{y}{x'}
= \qyx{y}{x}\qyx{y}{x'} \left(\nabla_{y}^{T} \log \qyx{y}{x}-\nabla_{y}^{T} \log \qyx{y}{x'} \right);
\end{equation}
\item for $x,x'$ such that $\qyx{y}{x}\qyx{y}{x'}=0$, we have $\qyx{y}{x}$ or (non-exclusive) $\qyx{y}{x'}=0$. For example assume $\qyx{y}{x}=0$. As $y' \mapsto \qyx{y'}{x}$ is non-negative, it is locally minimum at $y'=y$, and as it is differentiable this implies $\nabla_{y} \qyx{y}{x} = 0$. Similarly if $\qyx{y}{x'}=0$ then $\nabla_{y}\qyx{y}{x'} = 0$. As a result~\eqref{eq:Simplified} remains valid with the convention $\nabla_{y} \log q_{P} = 0$ where $q_{P}=0$.
\eit
With the above observations we rewrite
\begin{eqnarray}
q_{P}^2(y) Df_{P}(y)
&=&
\EXP_{X,X' }\Big(
X\qyx{y}{X} \qyx{y}{X'}
\left(
\nabla_{y}^{T} \log \qyx{y}{X} -  \nabla_{y}^{T} \log \qyx{y}{X'}
\right)\Big)\notag\\
&=&
\EXP_{X,X' } \Big(
-X\qyx{y}{X} \qyx{y}{X'}
\left(
\nabla_{y}^{T} \log \qyx{y}{X'} -  \nabla_{y}^{T} \log \qyx{y}{X}
\right)\Big)\label{eq:tmp1}\\
&\stackrel{X \leftrightarrow X'}{=}&
\EXP_{X,X'}  \Big(
-X'\qyx{y}{X'} \qyx{y}{X}
\left(
\nabla_{y}^{T} \log \qyx{y}{X} -  \nabla_{y}^{T} \log \qyx{y}{X'}
\right)\Big)\notag\\
&=&
\EXP_{X,X'} \Big(
X'\qyx{y}{X} \qyx{y}{X'}
\left(
\nabla_{y}^{T} \log \qyx{y}{X'} -  \nabla_{y}^{T} \log \qyx{y}{X}
\right)\Big)\label{eq:tmp2}
\end{eqnarray}
We conclude by taking the half sum of~\eqref{eq:tmp1} and~\eqref{eq:tmp2}.

\subsection{Proof of \prop{vaddmv}} \lab{pfvaddmv}
Since $p_{Y|X}$ is of the form \eq{expfam0} we have
$a(x)+b(y)=F(x-y)+c\<x,y\>$ for any $x,y \in \H$. A consequence is that $a(x)=F(x)-b(0)$ and $b(y)=F(-y)-a(0)$ for any $x,y$. In particular $a(0)+b(0)=F(0)$. It follows that $F(x-y)+c\<x,y\>=a(x)+b(y)=F(x)+F(-y)-a(0)-b(0)=F(x)+F(-y)-F(0)$ for any $x,y$.
Denoting $G(z) = G(z)-G(0)$, we get $G(x-y)+c\<x,y\>=G(x)+G(-y)$. Specializing to $x=y$ yields $G(0)+c\|x\|^2=G(x)+G(-x)$ for any $x$, hence $G(0)=0$ and $c \|x\|^2= G(x)+G(-x)$ for all $x$. Denoting $A(z),B(z)$ the odd and the even part of $G(z)$ we get $B(z) = c\|z\|^2/2$ and 
\[
A(x-y)+c\|x-y\|^2/2 + c\<x,y\> = G(x-y) + c \<x,y\> = G(x)+G(-y) = A(x)+A(-y)+c\|x\|^2/2+c\|y\|^2/2
\]
for any $x,y$. Thus, $A(x-y)=A(x)+A(-y)$ for any $x,y$ and, as $A$ is $C^0$ (by the continuity of $F$) it follows that $A$ is linear: there is $\mu \in \H$ such that $A(x) = -c\<x,\mu\>$ for all $x$, so that $G(z) = c\|z\|^2/2-c\<z,\mu\> = \tfrac{c}{2}\lp(\|z-\mu\|^2-\|\mu\|^2\rp)$ and $F(z) = c\|z-\mu\|^2+d$ with $d \in \RR$. As the noise is centered, we conclude that $\mu = 0$.

\subsection{Proof of \lem{cvxbnd}} \lab{pfcvxbnd}
Equivalently we show that $G(x)-\log(1+|x|)$ is bounded from below. 

Since $\int_\RR e^{-G(x)} dx < \infty$ and $G$ is convex, there is $a \in \RR$ such that $G$ is non-increasing on $(-\infty,a]$ and non-decreasing on $[a,+\infty)$ with $\lim_{|x| \to \infty} G(x) = +\infty$. By convexity again, it follows that there are $x_0 < a < x_1$ and $u_0 < 0 < u_1$ such that $G(x) \geq G(x_1)+u_1(x-x_1)$ and $G(x) \geq G(x_0)+u_0(x-x_0)$ for any $x \in \RR$. On $[x_1,+\infty)$ we have $G(x)-\log(1+|x|) \geq G(x_1)+u_1(x-x_1)-\log(1+|x|)$ hence
\[
\inf_{x \in [x_1,\infty)} \{G(x)-\log(1+|x|)\} \geq G(x_1)-u_1 x_1 + \inf_{x \geq x_1} \{u_1 x_1-\log(1+|x|)\} > -\infty. 
\]
Similarly $\inf_{x \in (-\infty,x_0]} \{G(x)-\log(1+|x|)\} > -\infty$. Finally, as $G$ is convex on $[x_0,x_1]$ it is continuous on this compact interval hence 
$\inf_{x \in [x_0,x_1]} \{G(x)-\log(1+|x|) = \min_{x \in [x_0,x_1]} \{G(x)-\log(1+|x|)\} > -\infty$. Putting the pieces together establishes the result.

\subsection{Worked example: scalar denoising with Laplacian noise and Laplacian prior} \lab{appscalL1L1}
Consider $p_{Y|X}(y|x) \propto \exp\lp(-|y-x|\rp) =: \qyx{y}{x}$ and $p_X(x) = c\exp\lp(-c|x|\rp)/2$.
Consider $y>0$ (the case $y<0$ is treated similarly by symmetry): we have
\begin{eqnarray*}
|y-x|+c|x| &=&
\begin{cases}
y-(1+c)x &\text{if}\ x < 0\\
y-(1-c)x &\text{if}\ 0 \leq x \leq y\\
(1+c)x-y &\text{if}\ x> y
\end{cases}
\end{eqnarray*}
hence for $c \neq 1$ we have
\begin{eqnarray*}
\tfrac{2}{c} q_P(y) &=& \int_{-\infty}^{+\infty} \exp\lp(-|y-x|-|x|\rp) dx\\
       &=& \int_{-\infty}^{0} e^{(1+c)x-y} dx 
       + \int_{0}^{y} e^{(1-c)x-y} dx
       + \int_{y}^{+\infty} e^{y-(1+c)x} dx\\
       &=&
       e^{-y} \int_{-\infty}^{0} e^{(1+c)x} dx
       + e^{-y} \int_{0}^{y} e^{(1-c)x}
       + \int_{0}^{+\infty} e^{y-(1+c)(y+x)} dx\\
   &=&
       e^{-y} \int_{0}^{+\infty} e^{-(1+c)x} dx
       + e^{-y} \int_{0}^{y} e^{(1-c)x}
       + e^{-cy} \int_{0}^{+\infty} e^{-(1+c)x} dx\\
       &=& \tfrac{e^{-y}+e^{-cy}}{1+c} 
       + e^{-y}\ \tfrac{e^{(1-c)y}-1}{1-c} 
       = 
       \tfrac{e^{-y}+e^{-cy}}{1+c} 
       + \tfrac{e^{-cy}-e^{-y}}{1-c}
       \end{eqnarray*}
        \begin{eqnarray*}
\tfrac{2}{c} \int_{-\infty}^{+\infty} x \qyx{y}{x} p_X(x)dx
       &=& \int_{-\infty}^{0} x\ e^{(1+c)x-y} dx 
       + \int_{0}^{y} x\ e^{(1-c)x-y} dx
       + \int_{y}^{+\infty} x\ e^{y-(1+c)x} dx\\
      &=& -e^{-y}\ \int_{0}^{+\infty} x\ e^{-(1+c)x} dx 
       + e^{-y}\ \int_{0}^{y} x\ e^{(1-c)x} dx\\
       &&+ \int_{0}^{+\infty} (y+x)\ e^{y-(1+c)(y+x)} dx\\
       &=& \lp(e^{-cy}-e^{-y}\rp)\ \int_{0}^{+\infty} x\ e^{-(1+c)x} dx 
       + e^{-y}\ \int_{0}^{y} x\ e^{(1-c)x} dx\\
       &&+ y e^{-cy}\ \int_{0}^{+\infty} e^{-(1+c)x} dx\\
       &=&
       \tfrac{e^{-cy}-e^{-y}}{(1+c)^2}
       +\tfrac{y e^{-y}}{1+c} 
       + e^{-y}\ \lp[\lp(x-\tfrac{1}{1-c}\rp)\ \tfrac{e^{(1-c)x}}{1-c}\rp]_{0}^{y}
       \\
       &=&
       \tfrac{e^{-cy}-e^{-y}}{(1+c)^2}
       +\tfrac{y e^{-y}}{1+c} 
       + \tfrac{\lp((1-c)y-1\rp)\ e^{-cy} + e^{-y}}{(1-c)^2}\\ 
\end{eqnarray*}
These derivations allow to express analytically $f(y) := \EXP(X|Y=y) = \frac{\int_{-\infty}^{+\infty} x \qyx{y}{x} p_X(x)dx}{q_P(y)}$. 

\bibliographystyle{plain}
\bibliography{biblio}

\end{document}

%% file: Var_Def.tex
\usepackage[latin1]{inputenc}

\usepackage{amsthm}
\usepackage{amsmath}
\usepackage{amssymb}
\usepackage{graphicx}
\usepackage{stackrel}
\usepackage{color}
\usepackage{enumerate}
\usepackage{epsfig}
\usepackage{algorithm}
\usepackage{booktabs,pstricks}
\usepackage{lipsum}
\usepackage{amsfonts}
\usepackage{graphicx}
\usepackage{epstopdf}
\usepackage{algorithmic}
\usepackage{amsopn}

\def\dst{\textcolor[rgb]{0.075,0.3,0.6059}}

\def\bit{\begin{itemize}}   \def\eit{\end{itemize}}
\def\ben{\begin{enumerate}} 
\def\een{\end{enumerate}}
\def\benA{\begin{enumerate}[\rm(A)]} 
\def\bena{\begin{enumerate}[\rm(a)]} 
\def\beni{\begin{enumerate}[\rm(i)]}

\def\beq{\begin{equation}}  \def\eeq{\end{equation}}
\def\beqn{\begin{eqnarray}} \def\eeqn{\end{eqnarray}}
\def\beqnn{\begin{eqnarray*}}   \def\eeqnn{\end{eqnarray*}}
\def\barr{\begin{array}}    \def\earr{\end{array}}
\def\bca{\begin{cases}}    \def\eca{\end{cases}}
\def\bfig{\begin{figure}}   \def\efig{\end{figure}}
\def\bpic{\begin{picture}}  \def\epic{\end{picture}}
\def\btab{\begin{tabular}}  \def\etab{\end{tabular}}

\def\ff{\theta}
\def\gg{\varrho}
\def\hh{h}

\def\bc{\begin{center}}     \def\ec{\end{center}}

\def\qu{\quad}

\def\ms{\medskip}

\def\RR{\mathbb R}
\def\NN{\mathbb N}

\def\EXP{\mathbb E}

\newtheorem{theorem}{\sc Theorem}
\newtheorem{proposition}{\sc Proposition}
\newtheorem{definition}{\sc Definition}
\newtheorem{lemma}{\sc Lemma}
\newtheorem{corollary}{\sc Corollary}
\newtheorem{remark}{\sc Remark}
\newtheorem{example}{\sc Example}
\newtheorem{fact}{\textsc{Fact}}

\def\BT{\begin{theorem}}      \def\ET{\end{theorem}}
\def\BP{\begin{proposition}}  \def\EP{\end{proposition}}
\def\BD{\begin{definition}}   \def\ED{\end{definition}}
\def\BL{\begin{lemma}}        \def\EL{\end{lemma}}
\def\BC{\begin{corollary}}    \def\EC{\end{corollary}}
\def\BR{\begin{remark}\rm }       \def\ER{\end{remark}} 
\def\BE{\begin{example}\rm }      \def\EE{\end{example}}
\def\BF{\begin{fact}}          \def\EF{\end{fact}}
\def\endproof{\hfill $\mb{\small$\Box$}$\ms}
\def\proof{\par\noindent{\textit{Proof}}. \ignorespaces}

\def\cor{Corollary~\ref}
\def\defi{Definition~\ref}
\def\rem{Remark~\ref}
\def\thm{Theorem~\ref}
\def\prop{Proposition~\ref}
\def\lem{Lemma~\ref}
\def\exa{Example~\ref}

\def\prox{\mathrm{prox}}
\newcommand\dom[1]{\operatorname{dom}(#1)}
\newcommand\Ima[1]{\operatorname{Im}(#1)}

\def\lab{\label}

\def\ed{\end{document}}
\def\eq{\eqref}
\def\lp{\left}  \def\rp{\right}
\def\<{\langle} \def\>{\rangle}
\def\leq{\leqslant} \def\geq{\geqslant}
\def\void{\varnothing}
\def\x{\times}
\def\all{\forall\;}

\def\#{\,\sharp\,}
\def\Id{\mathrm{Id}}

\def\dim{\mathrm{dim}\,}

\def\d{\partial}
\def\sign{\mathrm{sign}}

\def\wt{\widetilde}

\def\mb{\mbox}

   \def\H{{\mathcal H}}

   \def\S{\mathcal{S}}

 \def\V{\mathcal{V}} \def\X{\mathcal{X}}
\def\Y{\mathcal{Y}} \def\Z{\mathcal{Z}}

\def\ph{\varphi}

\newcommand\qyx[2]{q(#1\,|\,#2)}

%% file: part2.bbl
\begin{thebibliography}{10}

\bibitem{Amini:2013co}
Arash Amini, Ulugbek Kamilov, Emrah Bostan, and Michael~A Unser.
\newblock {Bayesian Estimation for Continuous-Time Sparse Stochastic
  Processes.}
\newblock {\em IEEE Trans. Signal Processing}, 61(4):907--920, 2013.

\bibitem{Amini:2011uh}
Arash Amini, Michael~A Unser, and Farokh Marvasti.
\newblock {Compressibility of Deterministic and Random Infinite Sequences}.
\newblock {\em IEEE Trans. Information Theory}, 59(11):5193--5201, 2011.

\bibitem{Banerjee:2005jd}
Arindam Banerjee, Xin~Guo 0001, and Hui~Wang 0003.
\newblock {On the optimality of conditional expectation as a Bregman
  predictor.}
\newblock {\em IEEE Trans. Information Theory}, 51(7):2664--2669, 2005.

\bibitem{Belge00}
Murat Belge, Misha Kilmer, and Eric Miller.
\newblock Wavelet domain image restoration with adaptive edge-preserving
  regularization.
\newblock {\em IEEE Trans. Image Process.}, 9(4):597--608, 2000.

\bibitem{Bregman67}
L.~M. Bregman.
\newblock The relaxation method of finding the common point of convex sets and
  its application to the solution of problems in convex programming.
\newblock {\em USSR Computational Mathematics and Mathematical Physics},
  7(3):200--217, 1967.

\bibitem{Burger:2014hv}
Martin Burger and Felix Lucka.
\newblock {Maximum a posteriori estimates in linear inverse problems with
  log-concave priors are proper Bayes estimators}.
\newblock {\em Inverse problems}, 30(11):114004--22, October 2014.

\bibitem{GRIBONVAL:2010:INRIA-00486840:1}
Remi Gribonval.
\newblock {Should Penalized Least Squares Regression be Interpreted as Maximum
  A Posteriori Estimation?}
\newblock {\em IEEE Transactions on Signal Processing}, 59(5):2405--2410, 2011.

\bibitem{gribonval:inria-00563207}
Remi Gribonval, Volkan Cevher, and Michael~E Davies.
\newblock {Compressible Distributions for High-Dimensional Statistics}.
\newblock {\em IEEE Trans. Information Theory}, 58(8):5016--5034, 2012.

\bibitem{NIPS2013_4868}
Remi Gribonval and Pierre Machart.
\newblock {Reconciling "priors" and "priors" without prejudice?}
\newblock In C~J~C Burges, L~Bottou, M~Welling, Z~Ghahramani, and K~Q
  Weinberger, editors, {\em Advances in Neural Information Processing Systems
  26 (NIPS)}, pages 2193--2201, 2013.

\bibitem{RGMN2018a}
R{\'e}mi Gribonval and Mila Nikolova.
\newblock {A characterization of proximity operators}.
\newblock {\em \rm\url{https://hal.inria.fr/hal-01835101}}, July 2018.

\bibitem{Helin:2015jq}
T~Helin and M~Burger.
\newblock {Maximum a posteriori probability estimates in infinite-dimensional
  Bayesian inverse problems}.
\newblock {\em Inverse problems}, 31(8):085009, August 2015.

\bibitem{SSP}
S~M Kay.
\newblock {\em {Fundamentals of Statistical Signal Processing{\textasciitilde}:
  Estimation Theory}}.
\newblock Signal Processing. Prentice Hall, 1993.

\bibitem{Kazerouni:2013co}
Abbas Kazerouni, Ulugbek Kamilov, Emrah Bostan, and Michael~A Unser.
\newblock {Bayesian Denoising - From MAP to MMSE Using Consistent Cycle
  Spinning.}
\newblock {\em IEEE Signal Process. Lett.}, 20(3):249--252, 2013.

\bibitem{Louchet:2013hs}
C{\'e}cile Louchet and Lionel Moisan.
\newblock {Posterior Expectation of the Total Variation Model: Properties and
  Experiments}.
\newblock {\em SIAM J. Imaging Sci.}, 6(4):2640--2684, January 2013.

\bibitem{Mathieu92}
P.~Mathieu, M.~Antonini, M.~Barlaud, and I.~Daubechies.
\newblock Image coding using wavelet transform.
\newblock {\em IEEE Trans. Image Process.}, 1(2):205--220, 1992.

\bibitem{Moreau65}
Jean-Jacques Moreau.
\newblock {Proximit{\'e} et dualit{\'e} dans un espace Hilbertien}.
\newblock {\em Bull. Soc. math. France}, 93:273--299, 1965.

\bibitem{Nikolova:2007aa}
Mila Nikolova.
\newblock {Model distortions in Bayesian MAP reconstruction}.
\newblock {\em Inverse Problems and Imaging}, 1(2):399--422, 2007.

\bibitem{Unser:2014vs}
Michael~A Unser and Pouya~D Tafti.
\newblock {\em {An introduction to sparse stochastic processes}}.
\newblock Cambridge University Press, 2014.

\end{thebibliography}
